\documentclass[10pt]{amsart}
\usepackage{amsfonts,amsmath,amssymb,amsthm,mathtools,enumerate,hyperref}
\usepackage[left=3cm,top=2cm,right=3cm,bottom=2cm]{geometry}
\allowdisplaybreaks

\theoremstyle{plain}
\newtheorem{theorem}{Theorem}[section]
\newtheorem{proposition}[theorem]{Proposition}
\newtheorem{lemma}[theorem]{Lemma}
\newtheorem{corollary}[theorem]{Corollary}
\newtheorem{conjecture}[theorem]{Conjecture}
\theoremstyle{definition}
\newtheorem{definition}[theorem]{Definition}
\newtheorem{remark}[theorem]{Remark}

\DeclareMathOperator{\Res}{Res}
\newcommand{\C}{\mathbb{C}}
\newcommand{\Z}{\mathbb{Z}}
\newcommand{\Ng}{N_{g}}
\newcommand{\eps}{\varepsilon}
\newcommand{\half}{\tfrac12}
\newcommand{\Sg}{\Sigma}
\newcommand{\Lg}{L_\eps}
\newcommand{\Ip}{\mathcal I_+}
\newcommand{\I}{\mathcal I}
\newcommand{\dd}{\partial}
\newcommand{\Gt}{\widetilde\Gamma}
\newcommand{\gt}{\widetilde\gamma}
\newcommand{\Nrm}[1]{\lVert #1\rVert}

\title[Szeg\H{o} quasiforms]
{Szeg\H{o} quasiforms:
a Schottky-scheme construction for vertex operator superalgebras} 

\author{A. Zuevsky}
\address{Institute of Mathematics, Czech Academy of Sciences, 
\v{Z}itn\'{a} 25, Prague, Czech Republic}
\email{zuevsky@yahoo.com}


\begin{document}

\begin{abstract}
We study the half weight ($N=\half$) analogue of the Bers quasiform on a Schottky group 
$\Gamma$ of genus $g$: the Poincar\'e series 
$S_\eps(x,y)=\sum_{\gt\in\Gt}j(\gt,x)/(\gamma x-y)\,(dx\,dy)^{1/2}$ over a lift 
$\Gt\subset\mathrm{SL}_2(\C)$ of $\Gamma$. The $\Ng:=2^{g}$ lifts give $\Ng$  
trivial $\alpha$-holonomy even spin structures $\Lg$. 
 The remaining $2^{2g}-\Ng$ spin structures
 are not treated in this paper.
  We prove, on an open set (H) of Schottky parameters and more generally 
whenever the critical exponent of $\Gamma$ is below $\half$ (any genus $g\ge 1$) the following:
 absolute and locally uniform convergence on $\Omega_0\times\Omega_0$, where $\Omega_0$ is the full region 
of discontinuity including the point $\infty$ (treated in homogeneous coordinates); 
  invariance without Eichler-type defect; skew-symmetry; a simple pole on the diagonal 
as the only singularity; $H^0(\Sg,\Lg)=0$; and hence identification of $S_\eps$ with the 
Szeg\H{o} kernel of $\Lg$. At genus one the kernel is computed explicitly and the relations 
between lifts and theta characteristics are determined. 
We further prove that the diagonal coefficient defines a projective connection that 
differs from the Bergman projective connection by a quadratic differential 
 constructed from $\theta[\delta](0)$, a residue-reconstruction formula for half-differentials 
with poles of any order, an expansion in the Schottky multipliers,  
and a formula for the free fermion stress tensor. 
\end{abstract}
\maketitle

\section{Introduction}

The Bers quasiform $\Psi_N(x,y)$ ($N\ge 2$; see \cite{TWb}) is a Poincar\'e series over 
a Schottky group built from $\pi_N=\frac{dx^Ndy^{1-N}}{x-y}$, multiplied by an auxiliary 
factor $\prod_{j=0}^{2N-2}\frac{y-A_j}{x-A_j}$.  The factor destroys  periodicity in 
$y$ and produces an Eichler-cocycle defect. The kernel has weight $(N,1-N)=(\half,\half)$ 
at $N=\half$, the weight of a Szeg\H{o} kernel, and there are $2N-1=0$ auxiliary points. 
This paper treats that case, without formal continuation in $N$.

\subsection{Relation to earlier work}

 Poincar\'e-type series $\sum_\gamma(\cdots)\gamma'^{\,w}$ over Schottky groups go back to 
Burnside and Poincar\'e. For integral weights $w\ge 1$ and their role in the theory of 
automorphic forms for Schottky groups see \cite{Be, Bu, Mk}. Schottky-series expressions 
are also standard in the operator formalism of string theory. The Bers quasiform on a 
Schottky-uniformized surface, its properties and use in expansions of meromorphic 
$N$-forms are developed in \cite{TWb}. This paper is the source of the genus $g$ Zhu recursion 
for vertex operator algebras of \cite{TW} (see Remark \ref{remTW} for what that recursion is). 
  The name ``Szeg\H{o} kernel'' also designates the reproducing 
kernel of Hardy spaces on domains in $\C^n$. For recent work in that setting  
(closed formulas and asymptotic expansions on balls and Reinhardt domains) see 
\cite{BlEn, MX}. Those are different objects from the spinor Szeg\H{o} kernel on a compact 
Riemann surface studied here (a section of $L\boxtimes L$ with a pole on the diagonal, 
not the reproducing kernel of a Hardy space), and we do not explore a relation. To our 
knowledge the following is specific to weight $\half$: the series needs no auxiliary points 
and is  invariant under both variables. It is a Szeg\H{o} kernel (the unique one, by 
$H^0=0$) for each of the $\Ng$ lifted spin structures. Its convergence, behaviour at 
$\infty$ and expansion in the Schottky multipliers can be controlled explicitly.

\subsection{Our results} 

Let us describe the results of this paper. 
We proved 
(modulo the quoted classical facts on theta functions and prime forms, \cite{Fa, Mu}): 
Sections 2-3 except Remark \ref{remdict}, Theorems \ref{thmred} and \ref{thmexp}, 
Proposition \ref{propT}. The estimates of Sections 3.1-3.2 and Theorem \ref{thmexp} 
are uniform on the domain (H); Theorem \ref{thmE} extends the qualitative results to every 
Schottky group satisfying the intrinsic condition (E) below.

We proved, depending on a hypothesis: Corollary \ref{corzhu} (on (A)) and Proposition 
\ref{proppf} (on (A$_{\rm ff}$)). 
 Neither hypothesis is established here for correlation 
functions obtained by a Schottky construction. These are statements about meromorphic 
sections of $\Lg$, not about the vertex operator superalgebra.

We conjecture: the dictionary between lifts and theta characteristics for 
$g\ge 2$ (see Remark \ref{remdict}); convergence and continuation beyond (E); Conjecture 
\ref{conjzhu} (general-VOSA genus $g$ recursion).

Not treated in this paper: the $2^{2g}-\Ng$ remaining spin structures 
(this set contains both even 
spin structures outside our family and all odd ones); zero modes; super-moduli.

\subsection{Notations}

 $s_\eps(x,y)$ is the scalar kernel in the coordinates $(x,y)$, 
$S_\eps=s_\eps(x,y)(dx\,dy)^{1/2}$ the corresponding section of $\Lg\boxtimes\Lg$, and 
$\hat s_\eps(\mathbf a,\mathbf b)$ ($\mathbf a$, $\mathbf b\in\C^2\setminus0$, Section 
\ref{secmain}) the bihomogeneous representative of the same section, which treats $\infty$ 
like any other point. Domains attached to fixed Schottky data $p=(w,\lambda)$ carry the 
subscript $p$. Statements uniform in the data are formulated on total spaces 
(see Section \ref{secspin}).

\subsection{Our hypotheses}

We assume the following:
 
 (H) $0<r<d/(8(2g-1))$: condition 
on the Schottky data, used for all quantitative estimates (Lemmas \ref{lemword}, 
\ref{lemchart}, Theorems \ref{thmconv}, \ref{thmexp}). 

(D): the discs $D_c$ are pairwise 
disjoint (so that $\Gamma$ is a classical Schottky group; used in Section 2).

 (E): $\sum_{\gt}\Nrm{\gt}^{-1}<\infty$, intrinsic. 

We have the following implications: (H) leads to (D) and 
(H) leads to (E) for every $g\ge 1$ (Theorem \ref{thmE} (iii)). Critical exponent  
$<\tfrac12$ leads to (E). 
  (E) together with (D) gives all qualitative results of Section 3, 
in every genus $g\ge 1$. Genus one does not require (H) since 
 Proposition \ref{propg1} treats it. 
 (A), (A$_{\rm ff}$) are hypotheses on vertex-algebraic correlators,  
used only in Section 4.2 and Section 6. They are not implied by, and do not imply, 
any of (H), (D), (E).

\subsection{Main Theorem} 

Let $\Gamma$ be a Schottky 
group of genus $g\ge 1$ with data $(w,\lambda)$ as in Section \ref{secspin}
satisfying (D) and (E) of Section \ref{secE} (for $g\ge 2$ this holds in particular 
under the open condition (H)), and let $\eps$ be one of the $\Ng$ lifts 
(see Proposition \ref{propcount}). Then:  
\emph{\begin{enumerate}[]
\item \textbf{Theorems \ref{thmconv}, \ref{thmE}}: 
$S_\eps(x,y)=\sum_{\gt\in\Gt}\pi(\gt x,y)$ converges absolutely and locally uniformly 
for $x$, $y\in\Omega_0\subset\widehat\C$ 
(the point $\infty$ and its orbit included, in homogeneous coordinates) 
with $y\notin\Gamma x$. Under (H) it is moreover jointly holomorphic in the parameters 
on the total space. The bound $s_\eps=O(1/|x|)$ as $x\to\infty$ holds only locally uniformly 
in $y\in\Omega_0\setminus\Gamma\infty$. 
\item \textbf{Propositions \ref{propequi}, \ref{propskew}:} 
$S_\eps(\gamma x,y)=j(\gt,x)^{-1}S_\eps(x,y)$, 
$S_\eps(x,\gamma y)=j(\gt,y)^{-1}S_\eps(x,y)$ 
and $S_\eps(y,x)=-S_\eps(x,y)$, with no additive Eichler-type defect. 
\item \textbf{Proposition \ref{propequi}:} $S_\eps$ descends to a global meromorphic section 
of $\Lg\boxtimes\Lg$ on $\Sg\times\Sg$. 
\item \textbf{Proposition \ref{propdiag} (ii):} its only singularity is a simple pole along 
the diagonal $\Delta$, with principal part the canonical section of 
$(\Lg\boxtimes\Lg)(\Delta)|_\Delta\cong\mathcal O_\Sg$. 
\item \textbf{Theorem \ref{thmnonsing}:} $H^0(\Sg,\Lg)=0$, and $S_\eps$ is the unique 
meromorphic section of $\Lg\boxtimes\Lg$ with properties (2)-(4). 
\item \textbf{Proposition \ref{propfay}, Corollary \ref{corberg}:} 
$S_\eps$ is the theta-function Szeg\H{o} kernel of the characteristic $\delta(\eps)$, 
and $R_\eps=12c_1$ is the Bergman projective connection 
 given by $6\sum_{i,j}\partial_i\partial_j\log\theta[\delta](0)\,\zeta_i\zeta_j$. 
\item \textbf{Proposition \ref{propg1}:} for $g=1$ the kernel and the relation 
$\eps\mapsto\delta(\eps)$ are given explicitly.
\end{enumerate}}

\begin{remark}
\label{remgenus}

Genus zero is empty since there is no Schottky group.  
 Sections 2-4 and Section 6 hold for every $g\ge 1$ under (D) and (E) 
(Theorem \ref{thmE}, Proposition \ref{propg1}). The quantitative part 
(Lemmas \ref{lemword}, \ref{lemchart}, Theorem \ref{thmconv} and Theorem \ref{thmexp}, 
whose constants use $r/d\le\frac1{24}$) is written for $g\ge 2$. At genus one 
the constants change, and Theorem \ref{thmE} and Proposition \ref{propg1} are used instead.
\end{remark}

The vertex operator superalgebra statements depend on the following assumptions: 
 the fermionic Zhu-type reduction and the free fermion Pfaffian formula 
are proved only under hypotheses (A), (A$_{\rm ff}$) on genus $g$ correlation 
functions which are not established for Schottky-sewn correlators, and 
the general VOSA genus $g$ recursion is a conjecture (see \cite{TW}). 

\section{Schottky data, lifts and spin structures}
\label{secspin}

Let $\I=\Ip\cup(-\Ip)$, $\Ip=\{1,\dots,g\}$. We deal with the following data: 
 points $w_c\in\C$ ($c\in\I$), and $\lambda_a\in\C^\times$ ($a\in\Ip$, 
 in particular $\lambda_a\ne 0$), 
with $\lambda_{-a}:=-\lambda_a$ and $\rho_a:=-\lambda_a^2$, 
$\rho_{-a}:=\rho_a$. Put $d:=\min_{c\ne c'}|w_c-w_{c'}|$, 
$r_c:=|\lambda_c|=|\rho_c|^{1/2}$, 
$r:=\max_c r_c$ (thus $|\rho|_{\max}=r^2$), and let (H)
 denote the condition

\begin{equation}
\tag{H}
0<r< \frac{d}{8(2g-1)}.
\end{equation}
The parameter space is $\{(w_c)_{c\in\I}\in\C^{2g}\ 
\text{distinct}\}\times(\C^\times)^{\Ip}$, 
and (H) singles out an open subset $\mathfrak P$ of it; 
  $\lambda_a=0$ is not a point of the parameter space (see Remark \ref{remlambda0}). 
 A point of $\mathfrak P$ is written as $p=(w,\lambda)$, and $d=d(p)$, $r_c=r_c(p)$, $r=r(p)$ 
are continuous functions on $\mathfrak P$. 
 
We have the following convention on parameter-dependent domains:
 every domain below built from $w_c$, $r_c$ 
 (namely $\mathcal K$, $\mathcal K_d$, $\overline{\mathcal D}$, $\Omega_0$) depends on $p$ 
and is written $\mathcal K_p$, etc., when this matters. 
  Statements claiming uniformity or
 joint holomorphy in $p$ are formulated on the total spaces, e.g., 
$\mathbb K:=\{(p,z)\in\mathfrak P\times\widehat\C:\ z=\infty\ 
\; \text{or}\; |z-w_c(p)|\ge r_c(p)/2\ \text{for all}\; c\}$, but never on a product  
$U\times\mathcal K\times\mathcal K$. 
 Let $\gamma_a(z)=w_{-a}+\rho_a/(z-w_a)$, $\gamma_{-a}:=\gamma_a^{-1}$, 
i.e., $\gamma_{-a}(u)=w_a+\rho_a/(u-w_{-a})$. The matrix 
\[
\gt_a:=\frac1{\lambda_a}\begin{pmatrix}w_{-a}&\rho_a-w_aw_{-a}\\1&-w_a\end{pmatrix}
\]
has determinant $-\rho_a/\lambda_a^2=1$, thus $\gt_a\in\mathrm{SL}_2(\C)$ lifts $\gamma_a$.
 For $\gt=\begin{pmatrix}a&b\\c&d\end{pmatrix}\in\mathrm{SL}_2(\C)$, 
$j(\gt,z):=(cz+d)^{-1}$ is, by definition, a factor of automorphy: 
a function on $\mathrm{SL}_2(\C)\times\C$ 
 (not a globally defined square root of $\gamma'$) 
 satisfying the cocycle identity $j(\gt\tilde\mu,z)=j(\gt,\mu z)j(\tilde\mu,z)$. 
 Only the induced action on sections (as in the spin bundle below) is globally meaningful. 
One has $j(\gt,z)^2=\gamma'(z)$ and $j(\gt^{-1},\gamma z)=j(\gt,z)^{-1}$, 
and for all $c\in\I$
\begin{equation}\label{eqjgen}
j(\gt_c,z)=\frac{\lambda_c}{z-w_c}.
\end{equation}
Thus,  the pair ($\rho_a$, sign $\eps_a$ relative to a fixed branch of $\rho_a^{1/2}$) 
is encoded in the single parameter $\lambda_a$. Changing $\eps_a$ is 
$\lambda_a\mapsto-\lambda_a$. Let $\Gt=\langle\gt_a\rangle\cong\Gamma$ 
(free, thus this is a homomorphic section; $\Gt$ and $\Gamma$ are kept distinct).

\subsection{Standard facts}

In this subsection we recall standard facts from \cite{Mk, Mu}.  
 Let $D_c=\{|z-w_c|\le r_c\}$ and consider the condition

\medskip\noindent(D)\quad the $2g$ discs $D_c$ are pairwise disjoint.

\medskip\noindent Condition (H) implies (D), since $d>2r$ there. All of this subsection 
and Proposition \ref{propcount} below use only (D). 

Let us assume (D). 
$\gamma_c$ maps $\{|z-w_c|>r_c\}$ into $D_{-c}$. Hence $\Gamma$ is a free classical 
Schottky group. $\overline{\mathcal D}:=\widehat{\C}\setminus\bigcup_cD_c^\circ$ is a 
closed fundamental domain with side pairings $\gamma_a:\partial D_a\to\partial D_{-a}$,
so that $\Omega_0=\Gamma\overline{\mathcal D}$.  
 Every point of the region of discontinuity $\Omega_0\subset\widehat\C$
 is a $\Gamma$-translate of a point of $\overline{\mathcal D}$. Boundary points
 may have two such representatives, which does not affect the existence statements used below.
  Note that $\overline{\mathcal D}\subset\widehat\C$ contains $\infty$. 
Hence, $\infty\in\Omega_0$, and the orbit $\Gamma\infty\subset\Omega_0$ consists 
of $\infty$ together with the finite points $\gamma(\infty)$, $\gamma\ne 1$  
(for instance $\gamma_c(\infty)=w_{-c}$). The image of $\infty$ in $\Sg$ is an ordinary 
point of $\Sg$. The scalar trivialization by $(dz)^{1/2}$ used below is singular 
along $\Gamma\infty$, which is why $\infty$ is treated in homogeneous coordinates 
in Section \ref{secmain}. $\Sg=\Omega_0/\Gamma$ has genus $g$.
 Every $\gamma\ne 1$ in $\Gamma$ is loxodromic with both fixed points in the limit set 
$\Lambda$, thus the action of $\Gamma$ on $\Omega_0$ is free and properly discontinuous 
(standard for classical Schottky groups, loc. cit.). 
Consequently, for fixed $x\in\Omega_0$, the orbit $\Gamma x$ is a locally finite subset 
of $\Omega_0$. It accumulates only on the limit set $\Lambda\subset\bigcup_cD_c^\circ$, 
which is disjoint from $\Omega_0$. Each point of $\Gamma x$ has a unique preimage 
group element, since the action is free. We will also use proper discontinuity 
in the following form: 
\emph{for compact $K$, $K'\subset\Omega_0$ the set 
$\{\gamma\in\Gamma:\gamma K\cap K'\ne\emptyset\}$ is finite.} 
In particular, the family of graphs $\{(x,\gamma x)\}_{\gamma\in\Gamma}$ 
in $\Omega_0\times\Omega_0$ is locally finite, not only each single orbit. 
Every point of $\Omega_0\times\Omega_0$ has a neighbourhood $K\times K'$  
meeting only finitely many graphs.

\subsection{The spin bundle}

 $\Lg$ is the line bundle on $\Sg$ whose sections are $\psi(z)(dz)^{1/2}$ 
on $\Omega_0\setminus\Gamma\infty$ with $\psi(\gamma z)j(\gt,z)=\psi(z)$ for
 $\gt\in\Gt$. The trivialization by $(dz)^{1/2}$ is singular at $\Gamma\infty$. 
 Across these points a section is extended by the transformation law, equivalently, 
 (see Section \ref{secmain})  $\Lg$ is the quotient by $\Gt\subset\mathrm{SL}_2(\C)$ of 
the restriction to $\Omega_0$ of the $\mathrm{SL}_2(\C)$-equivariant bundle 
$K^{1/2}_{\mathbb P^1}=\mathcal O(-1)$. Since $j^2=\gamma'$,
 $\Lg^{\otimes2}\cong K_\Sg$, and $\deg\Lg=g-1$. 
The next proposition concerns the family of spin structures. 

\begin{proposition}
\label{propcount}
Let us assume the condition (D). Then: 

 \emph{(a)} The $\Ng$ lifts of $\Gamma$ (equivalently, the $\Ng$ 
sign choices $\lambda_a\mapsto\pm\lambda_a$) give $\Ng$ pairwise non-isomorphic 
spin structures $\Lg$. 

\emph{(b)} They are the spin structures $L$ for 
which the pull-back $\pi^*L$ to $\Omega_0$ is the trivial spin structure, i.e., 
for which $dz$ has a single valued square root $(dz)^{1/2}$ as a section of $\pi^*L$ over 
$\Omega_0$.

\emph{(c)} For the marking with $\alpha_a=\partial D_a$ their characteristics are 
$\delta=[0;\delta'']$, $\delta''\in\{0,\half\}^g$, all even
(proven in Proposition \ref{propfay} (b)). 
\end{proposition} 

We do not discuss here which characteristic corresponds to which sign vector for $g\ge 2$  
(see Remark \ref{remdict}). The remaining $2^{2g}-\Ng$ spin structures,  
both even structures not of this form and every odd structure, are not treated. 
 Note that in \emph{(b)} above we do not claim that $\pi^*L$ is trivial as a 
holomorphic line bundle, 
which holds for every $L$ because $\Omega_0$ is a plane domain, nor does every spin  
structure admit such a square root.

\begin{proof}
(a) \emph{Lifts:} a lift of $\gamma_a$ to $\mathrm{SL}_2(\C)$ is $\pm\gt_a$, 
i.e., $\lambda_a\mapsto\pm\lambda_a$. Since $\Gamma$ is free on 
$\gamma_1$, $\dots$, $\gamma_g$, every assignment of lifts of the generators extends 
uniquely to a homomorphic section $\Gamma\to\mathrm{SL}_2(\C)$. 
 This gives the $\Ng$ groups $\Gt$. 

{\it Bundles:} since $j(\gt,z)^2=\gamma'(z)$ and $j$ satisfies the cocycle identity, 
$\gt\cdot(z,v):=(\gamma z,j(\gt,z)^{-1}v)$ is an action on $\Omega_0\times\C$, free and 
properly discontinuous on the base. $\Lg:=(\Omega_0\times\C)/\Gt$ is a line bundle 
on $\Sg$ with $\Lg^{\otimes2}\cong K_\Sg$. The section $v=1$ trivializes the pull-back, 
and corresponds to $(dz)^{1/2}$. Two lifts differ by a character
 $\chi:\Gamma\to\{\pm1\}$, and $L_{\eps'}\cong L_\eps\otimes F_\chi$ 
with $F_\chi$ the flat unitary line bundle of $\chi$, viewed as a character 
of $\pi_1(\Sg)$ through the surjection $\pi_1(\Sg)\to\Gamma$. 
Distinct unitary characters of $\pi_1(\Sg)$ give non-isomorphic holomorphic line bundles, 
 (injectivity of the map from unitary characters to $\mathrm{Pic}^0$),
thus the $\Ng$ bundles $\Lg$ are pairwise non-isomorphic. 

(b) Spin structures on $\Sg$ form a torsor under $H^1(\Sg;\Z/2)$, 
and the restriction of a spin structure to $\Omega_0$ is trivial iff its class 
pulled back to $\Omega_0$ coincides with that of $(dz)^{1/2}$. 
 I.e., iff its holonomy along each generator $\alpha_a=\partial D_a$ of
 $\pi_1(\Omega_0)$ equals that of $(dz)^{1/2}$ 
(which is trivial, $(dz)^{1/2}$ being single valued on $\Omega_0$). 
Namely, $2^{2g}/2^g=\Ng$ classes satisfy this ($g$ independent conditions),
 and similarly the $\Lg$.

  (c) is Proposition \ref{propfay} (b).
\end{proof}

\subsection{Convention} 

For the fixed marking ($\alpha_a=\partial D_a$,  
$\beta_a$ as in Proposition \ref{propfay}), theta characteristics 
$\delta\in(\tfrac12\Z/\Z)^{2g}$ correspond bijectively to spin structures $L_\delta$ 
(the spin structure whose Szeg\H{o} kernel, when $\theta[\delta](0)\ne 0$, 
is given by the formula of Proposition \ref{propfay}). For each lift $\eps$ 
we let $\delta(\eps)$ be the unique characteristic with $L_{\delta(\eps)}\cong\Lg$. 
This is defined purely through the isomorphism class of the bundle $\Lg$ 
(i.e., no choice of sign or branch). 
  In what follows, we use $\delta(\eps)$ of this form.  
The numerical relation between the sign vector of $\eps$ and $\delta''(\eps)$ is 
a separate question (Proposition \ref{propg1} for $g=1$, Remark \ref{remdict} for $g\ge 2$) 
on which no proof in this paper depends. We also mention the classical fact that 
(Riemann singularity theorem in its form for characteristics, see \cite{Mu, Fa})
\[
\theta[\delta](0)=0\;\;  \textit{iff} \;\;  h^0(\Sg,L_\delta)>0,
\]
valid for every characteristic $\delta$ (even or odd), and  we do prove it again here. 

\section{The Szeg\H{o} quasiform}
\label{secmain}

Let $\pi(x,y)={(dx\,dy)^{1/2}}(x-y)^{-1}$. 
 By $gx-gy=\frac{x-y}{(cx+d)(cy+d)}$ it is jointly invariant, 
$\pi(\tilde gx,\tilde gy)=\pi(x,y)$ for $\tilde g\in\mathrm{SL}_2(\C)$. 
In coordinates: $j(\tilde g,x)j(\tilde g,y)/(gx-gy)=1/(x-y)$.

\begin{definition}
\label{defS}
For $x$, $y\in\Omega_0$ with $y\notin\Gamma x$, 
\[
s_\eps(x,y):=\sum_{\gt\in\Gt}\frac{j(\gt,x)}{\gamma x-y},
\qquad S_\eps:=s_\eps(x,y)(dx\,dy)^{1/2}. 
\]
\end{definition}

We will use the following conventions: $\gamma$ denotes the M\"obius transformation 
underlying $\gt$, $j(\gt,x)=(c_{\gt}x+d_{\gt})^{-1}$, and the first variable $x$ is
 the one acted on by $\Gt$, the second variable $y$ the one at which the pole sits 
($y=\gamma x$). The homogeneous formulation below (eq. \eqref{eqshat}) is the same object. 
 With $\mathbf x=(x,1)$, $\mathbf y=(y,1)$ one has  
$\langle\gt\mathbf x,\mathbf y\rangle=j(\gt,x)^{-1}(\gamma x-y)$, 
thus the term $1/\langle\gt\mathbf x,\mathbf y\rangle$ equals 
$j(\gt,x)/(\gamma x-y)$ with no sign. The pairing is skew symmetric,   
$\langle\mathbf x,\mathbf y\rangle=-\langle\mathbf y,\mathbf x\rangle$, 
which is the source of the skew-symmetry of Proposition \ref{propskew} (i). 
 A change of chart $(x,1)=x\,(1,u)$, $u=1/x$, is absorbed by bihomogeneity 
of degree $(-1,-1)$ (the constant $\pm i$ coming from $dx=-u^{-2}du$ enters 
only the identification of $K^{1/2}$ with $\mathcal O(-1)$ and is the same for all terms). 
That this series converges on the stated set is the content of Theorem \ref{thmconv} 
(until then it is only a formal expression). The definition makes sense only for finite 
$x$, $y$, and the point $\infty\in\Omega_0$ is dealt with by the following reformulation.

\subsection{Homogeneous coordinates (treatment of $\infty$)} 

For $\mathbf x=(x_1,x_2)$, $\mathbf y=(y_1,y_2)\in\C^2\setminus0$ 
set $\langle\mathbf x,\mathbf y\rangle:=x_1y_2-x_2y_1$. 
 Then $\langle\tilde g\mathbf x,\tilde g\mathbf y\rangle=\langle\mathbf x,\mathbf y\rangle$ 
for $\tilde g\in\mathrm{SL}_2(\C)$. The point $x\in\C$ is represented by $\mathbf x=(x,1)$, 
the point $\infty$ by $\mathbf x=(1,0)$, and a point $\mathbf x$ of the chart at infinity by 
$(1,u)$, $u=1/x$. A section of $K^{1/2}_{\mathbb P^1}\cong\mathcal O(-1)$ 
over $V\subset\widehat\C$ is a function $f$ on the cone over $V$ 
with $f(t\mathbf x)=t^{-1}f(\mathbf x)$. In the chart $\mathbf x=(x,1)$ it is 
$f(x,1)(dx)^{1/2}$ and in the chart $\mathbf x=(1,u)$ it is $f(1,u)(du)^{1/2}$ 
(a constant factor $\pm i$ from $dx=-u^{-2}du$ is absorbed in the identification 
$K^{1/2}\cong\mathcal O(-1)$ and plays no role). For 
$\tilde g=\begin{pmatrix}\alpha&\beta\\\kappa&\varsigma\end{pmatrix}\in\mathrm{SL}_2(\C)$ 
one has $\tilde g(x,1)=(\kappa x+\varsigma)\,(gx,1)=j(\tilde g,x)^{-1}(gx,1)$. 
 Hence $f(\tilde g\mathbf x)=j(\tilde g,x)\,f(gx,1)$ for $\mathbf x=(x,1)$.
 Invariance $f(\tilde g\mathbf x)=f(\mathbf x)$ is the automorphy law 
$\psi(gx)j(\tilde g,x)=\psi(x)$ defining $\Lg$. 
 Thus, $\Lg=\mathcal O(-1)|_{\Omega_0}/\Gt$, and sections of $\Lg$ are 
the $\Gt$-invariant functions $f$. Define, for $\mathbf x$, $\mathbf y\in\C^2\setminus0$
 over $\Omega_0$ with $[\mathbf y]\notin\Gamma[\mathbf x]$,

\begin{equation}
\label{eqshat}
\hat s_\eps(\mathbf x,\mathbf y)
:=\sum_{\gt\in\Gt}\frac1{\langle\gt\mathbf x,\mathbf y\rangle},
 \qquad \hat s_\eps(t\mathbf x,t'\mathbf y)=t^{-1}t'^{-1}\hat s_\eps(\mathbf x,\mathbf y).
\end{equation}
For $\mathbf x=(x,1)$, $\mathbf y=(y,1)$ one computes 
 $\langle\gt\mathbf x,\mathbf y\rangle=j(\gt,x)^{-1}(\gamma x-y)$. Thus, 
$\hat s_\eps((x,1),(y,1))=s_\eps(x,y)$ is the series of Definition \ref{defS}. 
 At $\mathbf x=(1,0)$ (the point $\infty$) $\hat s_\eps$ is the same series 
in the other chart, and 
$\langle\gt(1,0),(y,1)\rangle=\mathrm{(first\ column\ of\ }\gt)\ \text{paired with }(y,1)
=a_{\gt}-c_{\gt}y=c_{\gt}(\gamma\infty-y)$ 
(for $\gt=\begin{pmatrix}a_{\gt}&\ast\\c_{\gt}&\ast\end{pmatrix}$, 
$c_{\gt}\ne 0$ for $\gt\ne 1$).

Let $\mathcal K=\mathcal K_p=\{z\in\C:|z-w_c|\ge r_c/2\ \textit{for all}\; c\}$, 
 $\widehat{\mathcal K}:=\mathcal K\cup\{\infty\}\subset\widehat\C$, a compact set containing 
$\overline{\mathcal D}$ in its interior, and $\varpi:=2r/d\le\frac1{4(2g-1)}$. 
 Set $W:=\max_c|w_c|$ and $X_0:=2(W+r)$.  
 A \emph{normalized representative} of a point $\alpha\in\widehat{\mathcal K}$ 
is $\mathbf x=(x,1)$ if $\alpha=x$ is finite with $|x|\le 2X_0$, 
and $\mathbf x=(1,u)$ with $|u|\le 1/X_0$ ($u=0$ for $\alpha=\infty$) 
if $|\alpha|\ge X_0$ or $\alpha=\infty$. 
 Points with $X_0\le|\alpha|\le 2X_0$ have both. 
 Two charts are then needed: no single normalization of $\mathbf x$ is bounded 
near $\infty$ and at the same time bounded away from $0$ on all of $\widehat{\mathcal K}$ 
using a single formula. 

Let $x\in\mathcal K$ and let $\gt=\gt_{c_1}\cdots\gt_{c_n}\in\Gt$ ($n\ge 1$) be a reduced word, 
i.e., $\gamma=\gamma_{c_1}\circ\cdots\circ\gamma_{c_n}$ with $c_{i+1}\ne -c_i$ for 
$1\le i\le n-1$ (so that $\gamma$ first applies $\gamma_{c_n}$ to $x$, 
then $\gamma_{c_{n-1}}$, and so on, ending with $\gamma_{c_1}$). Then we have 

\begin{lemma}
\label{lemword}
 $|j(\gt,x)|\le 2\varpi^{n-1}$. Moreover $|\gamma x-w_{-c_1}|\le 2r_{c_1}$ for $n=1$ and 
$\le r_{c_1}\varpi\le r_{c_1}/4$ for $n\ge 2$. There are $2g(2g-1)^{n-1}$ reduced words 
of length $n$.
\end{lemma}

\begin{proof}
Set $z_n:=x$ and, for $i=n,n-1$, $\dots$, $1$, $z_{i-1}:=\gamma_{c_i}(z_i)$.  
Thus, $z_0=\gamma_{c_1}(\gamma_{c_2}(\cdots\gamma_{c_n}(x)\cdots))=\gamma x$. 
 By the cocycle identity $j(\gt,z_n)=\prod_{i=1}^nj(\gt_{c_i},z_i)$ 
 (the factor $j(\gt_{c_i},z_i)$ is evaluated at $z_i$, the point on which 
$\gt_{c_i}$ acts before $\gt_{c_{i-1}}$, $\dots$, $\gt_{c_1}$ are applied). 
 We estimate the factors in the order $i=n$, $n-1$, $\dots$, $1$. 
  I.e., starting from $x=z_n$ and working outward. 
 By \eqref{eqjgen}, $|j(\gt_{c_n},z_n)|=r_{c_n}/|z_n-w_{c_n}|\le 2$ 
 (as $z_n=x\in\mathcal K$, $|z_n-w_{c_n}|\ge r_{c_n}/2$) 
and $|z_{n-1}-w_{-c_n}|=r_{c_n}^2/|z_n-w_{c_n}|\le 2r_{c_n}$. 
 For $1\le i\le n-1$: since $c_i\ne -c_{i+1}$ and $|w_{-c_{i+1}}-w_{c_i}|\ge d$, 
 the bound $|z_i-w_{-c_{i+1}}|\le 2r_{c_{i+1}}\le 2r$ already established gives 
 $|z_i-w_{c_i}|\ge d-2r\ge d/2$ by (H). 
Then $|j(\gt_{c_i},z_i)|=r_{c_i}/|z_i-w_{c_i}|\le 2r_{c_i}/d\le\varpi$ and 
$|z_{i-1}-w_{-c_i}|=r_{c_i}^2/|z_i-w_{c_i}|\le r_{c_i}\varpi\le r_{c_i}/4$ 
(using $\varpi\le\frac1{4(2g-1)}\le\frac14$). 
Multiplying the $n$ bounds $|j(\gt_{c_n},z_n)|\le 2$ and 
$|j(\gt_{c_i},z_i)|\le\varpi$ ($i<n$) gives $|j(\gt,x)|\le 2\varpi^{n-1}$. 
  The case $i=1$ of the second estimate gives 
$|z_0-w_{-c_1}|=|\gamma x-w_{-c_1}|\le r_{c_1}\varpi\le r_{c_1}/4$ 
for $n\ge 2$ (for $n=1$,  $|\gamma x-w_{-c_1}|\le 2r_{c_1}$). 
Counting reduced words of length $n$ over the alphabet $\I$ ($|\I|=2g$) 
with no two consecutive inverse letters is the standard count $2g(2g-1)^{n-1}$.
\end{proof}

The next lemma describes the chart-normalized word estimate. 

\begin{lemma}
\label{lemchart}
Let $\alpha$, $\beta\in\widehat{\mathcal K}$ with normalized representatives 
$\mathbf x$, $\mathbf y$, and let $\gt=\gt_{c_1}$ $\cdots$ $\gt_{c_n}\in\Gt$ 
be a reduced word of length $n\ge 2$. Then
\[
\bigl|\langle\gt\mathbf x,\mathbf y\rangle\bigr|^{-1}\ \le\ C_*\,\varpi^{\,n-2},
\qquad C_*:=\frac{16}{d}+1+r.
\]
\end{lemma}

Note that the constant depends on the data only through $d$ and $r$. 
It does not involve $r_{\min}$, 
$\alpha$ or $\beta$. It is not coordinate-invariant. 
 The normalized representatives $(x,1)$, 
$(1,u)$ refer to the fixed affine coordinate on $\widehat\C$. Thus, $C_*$ adds quantities 
of different dimensions, $1/d$ and $1+r$. But this is not a problem since it  
is repaired by measuring  lengths in units of $d$. 

\begin{proof}

Set $z_n:=\alpha$ and $z_{i-1}:=\gamma_{c_i}(z_i)$ for $i=n$, $\dots$, $1$. 
Since $\gamma_c$ maps $\widehat{\C}\setminus D_c^\circ$ into $D_{-c}\subset\C$, 
all $z_i$ with $i\le n-1$ are finite. For $\mathbf v=(v_1,v_2)\in\C^2$ and $c\in\I$,
\[
\gt_c\mathbf v=\frac1{\lambda_c}\Bigl(w_{-c}v_1+(\rho_c-w_cw_{-c})v_2,\ v_1-w_cv_2\Bigr),
\]
so that $\gt_c\bigl(\mu\,(z,1)\bigr)=\mu\,\frac{z-w_c}{\lambda_c}\,(\gamma_cz,1)$. 
Applying this $n$ times (the first time to $\mathbf x$),

\begin{equation}
\label{eqchain}
\gt\mathbf x=\mu_1\mu_2\cdots\mu_n\,(z_0,1),\quad 
\mu_i=\frac{z_i-w_{c_i}}{\lambda_{c_i}}, \quad 1\le i\le n-1, 
\end{equation}
where $\mu_n$ is the second component of $\gt_{c_n}\mathbf x$, $z_0=\gamma x$. 

\emph{(a) The first part.} If $\mathbf x=(x,1)$ with $x\in\mathcal K$, 
then $\mu_n=(x-w_{c_n})/\lambda_{c_n}$. 
Thus, $|\mu_n|\ge\frac12$, and $|z_{n-1}-w_{-c_n}|=r_{c_n}^2/|x-w_{c_n}|\le 2r_{c_n}$. 
If $\mathbf x=(1,u)$ with $|u|\le 1/X_0$ (then $\alpha=1/u$ belongs outside all $D_c$, 
 since $X_0>W+r$), 
 then $\mu_n=(1-uw_{c_n})/\lambda_{c_n}$ with $|1-uw_{c_n}|\ge 1-W/X_0\ge\frac12$. 
 Thus $|\mu_n|\ge 1/(2r_{c_n})\ge 1/(2r)$. 
Moreover, $z_{n-1}-w_{-c_n}=\rho_{c_n}u/(1-uw_{c_n})$. 
Hence, $|z_{n-1}-w_{-c_n}|\le 2r_{c_n}^2/X_0\le r_{c_n}$ as $X_0\ge 2 r$.

\emph{(b) The later parts.} Let $1\le i\le n-1$ and assume $|z_i-w_{-c_{i+1}}|\le 2r$ 
 which is true for $i=n-1$ by (a). 
 Since $c_i\ne -c_{i+1}$ we have $|w_{c_i}-w_{-c_{i+1}}|\ge d$. 
 Thus, $|z_i-w_{c_i}|\ge d-2r\ge d/2$ by (H). 
 Hence, $|\mu_i|\ge d/(2r_{c_i})\ge d/(2r)=\varpi^{-1}$ and 
\[
|z_{i-1}-w_{-c_i}|=\frac{r_{c_i}^2}{|z_i-w_{c_i}|}\le\frac{2r_{c_i}^2}{d}
\le r_{c_i}\varpi\le\frac{r_{c_i}}4\ \ (\le 2r),
\]
which is the hypothesis for $i-1$. In particular, $|z_0-w_{-c_1}|\le r_{c_1}/4$ and 
$|\mu_1|\ge d/(2r_{c_1})$.

\emph{(c) The pairing with $\mathbf y$.} If $\mathbf y=(y,1)$ with $y\in\mathcal K$, 
then $\langle(z_0,1),\mathbf y\rangle=z_0-y$, and $|y-w_{-c_1}|\ge r_{c_1}/2$ together 
with (b) gives $|z_0-y|\ge r_{c_1}/4$. If $\mathbf y=(1,v)$ with $|v|\le 1/X_0$, 
then $\langle(z_0,1),\mathbf y\rangle=z_0v-1$ and $|z_0|\le W+r=X_0/2$, 
thus $|z_0v-1|\ge\frac12$.

\emph{(d) Conclusion.} By \eqref{eqchain}, 
$|\langle\gt\mathbf x,\mathbf y\rangle|
=|\mu_n|\cdot\prod_{i=2}^{n-1}|\mu_i|\cdot|\mu_1|\,|\langle(z_0,1),\mathbf y\rangle|$, 
 where the middle product has $n-2$ factors, each greater or equal $\varpi^{-1}$. 
By (b),(c): 
$|\mu_1||\langle(z_0,1),\mathbf y\rangle|\ge\frac d{2r_{c_1}}\cdot\frac{r_{c_1}}4=\frac d8$ 
if $\mathbf y=(y,1)$, and $\ge\frac d{2r}\cdot\frac12=\frac d{4r}$ if $\mathbf y=(1,v)$. 
Together with (a) this gives, for the four combinations of charts,
\[
|\langle\gt\mathbf x, \mathbf y\rangle|^{-1}\le\varpi^{n-2}\times
\begin{cases}2\cdot\frac8d=\frac{16}d,&\mathbf x, \; \mathbf y\ \text{both in the first chart},
\\ 2\cdot\frac{4r}d\le\frac13,&\mathbf x\ \text{first},\ \mathbf y\ \text{second},
\\ 2r\cdot\frac8d\le\frac23,&\mathbf x\ \text{second},\ \mathbf y\ \text{first},
\\ 2r\cdot\frac{4r}d\le\frac r3,&\text{both second},\end{cases}
\]
where we used $r/d\le\frac1{24}$ from (H). 
Each right-hand side is little or equal $C_*\varpi^{n-2}$. 
The four cases are collected in the following table, where $L_{\mathbf x}$ is the 
lower bound for $|\mu_n|$ from (a), $L_{\mathbf y}$ the lower bound for 
$|\mu_1|\,|\langle(z_0,1),\mathbf y\rangle|$ from (b)-(c), and the bound for 
$|\langle\gt\mathbf x,\mathbf y\rangle|^{-1}/\varpi^{n-2}$ is 
$1/(L_{\mathbf x}L_{\mathbf y})$:
\[
\begin{array}{c|c|c|c|c}
\mathbf x & \mathbf y & L_{\mathbf x}\ (\text{from (a)}) & L_{\mathbf y}\ 
(\text{from (b),(c)}) & 1/(L_{\mathbf x}L_{\mathbf y})\\\hline
(x,1) & (y,1) & 1/2 &  d/8 & {16}/d\\
(x,1) & (1,v) & 1/2 & d/{4r} & {8r}/d\le 1/3\\
(1,u) & (y,1) & 1/{2r} & d/8 & {16r}/d\le 2/3\\
(1,u) & (1,v) & 1/{2r} & d/{4r} & {8r^2}/d\le r/3 
\end{array}
\]
The second and fourth rows use $|\mu_1|\ge d/(2r)$ and $|z_0v-1|\ge\frac12$. 
The first and third use $|\mu_1|\ge d/(2r_{c_1})$ and $|z_0-y|\ge r_{c_1}/4$.
\end{proof}

\begin{remark}
\label{remconsts}

In this Remark we discuss constants and conventions in Lemmas 
\ref{lemword}-\ref{lemchart}. 
 The bounds $|j(\gt,x)|\le 2\varpi^{n-1}$ and 
$|\langle\gt\mathbf x,\mathbf y\rangle|^{-1}\le C_*\varpi^{n-2}$ hold for every reduced word 
of length $n$ with constants depending only on the Schottky data through $d$, $r$ and $\varpi$
 ($C_*=16/d+1+r$). They do not depend on the compact subset or on the word, and no constant
 has been enlarged to absorb short words. The identity element ($n=0$) is excluded from the 
reduced words of the lemmas and is treated separately. It contributes the term 
$1/\langle\mathbf x,\mathbf y\rangle=1/(x-y)$ (which carries the diagonal pole), and the words 
of length $1$ are likewise summed separately (see Lemma \ref{lemword} with $n=1$). 
In Theorem \ref{thmconv} the tail $\sum_{n\ge 2}$ is thus dominated by 
$C_*\sum_{n\ge 2}2g(2g-1)^{n-1}\varpi^{n-2}\le C_*\,2g(2g-1)\sum_{k\ge 0}q^k$, 
$q=(2g-1)\varpi\le\frac14$. The only dependence on the compact set $\widehat{\mathcal K}_p$ 
is through the finitely many terms $n\le 1$ and the distance to the polar loci. 
The same convention is used in Theorem \ref{thmexp}, whose remainder constant 
$22g(2g-1)|\rho|_{\max}/d^3$ is uniform on $\mathbb G_d$. The compact set enters only 
through $\mathcal K_{d,p}$, i.e., distance $\ge d/2$ from all $w_c$. It is not claimed 
it is uniform up to the pole loci $x\in\Gamma y$ (the expansion is stated on 
$\mathcal K_{d,p}\times\mathcal K_{d,p}$ only).
\end{remark}

The theorem below states the convergence and meromorphicity on $\Omega_0\times\Omega_0$, 
including $\infty$. Let us assume (H). Let $c_{\gt}$ denote the lower-left entry 
of $\gt\in\Gt$ 
(thus $c_{\gt}\ne 0$ for $\gt\ne 1$, since $\Gamma$ acts freely and $\gamma\infty\ne\infty$).

\begin{theorem}
\label{thmconv}
 
\begin{enumerate}[(i)]

\item {\rm Uniform tail on the total space.} Let $U\subset\mathfrak P$ be compact. 
For every $p\in U$ and all $\alpha$, $\beta\in\widehat{\mathcal K}_p$ with normalized 
representatives $\mathbf x$, $\mathbf y$,
\[
\sum_{n\ge 2}\ \sum_{|\gt|=n}\bigl|\langle\gt\mathbf x,\mathbf y\rangle\bigr|^{-1}\ \le\ 
\tfrac83\,g(2g-1)\,\sup_{U}C_*\ <\infty.
\]
Let $\mathbb U$ denote the interior of the total space 
$\{(p,\alpha,\beta):p\in\mathfrak P,\ \alpha,\beta\in\widehat{\mathcal K}_p\}$, 
in the coordinates $(x,y)$, $(u,y)$, $(x,v)$ or $(u,v)$, $u=1/x$, $v=1/y$, 
according to the charts. 
 Let $\mathbb P_1:=\{\beta=\alpha\}$, $\mathbb P_c:=\{\beta=\gamma_c\alpha\}$, $c\in\I$, 
be the polar graphs in $\mathbb U$. Then:

\begin{enumerate}[(a)]

\item the tail $T:=\sum_{|\gt|\ge 2}\langle\gt\mathbf x,\mathbf y\rangle^{-1}$ converges 
absolutely and uniformly on $\mathbb U$ (on the part over $U$) and is holomorphic on all 
of $\mathbb U$, including across the polar graphs (Weierstrass);

\item the terms $\gt=1$ and $\gt=\gt_c$ ($c\in\I$) are explicit rational functions of 
$(p,\mathbf x,\mathbf y)$, holomorphic on $\mathbb U\setminus\bigcup\mathbb P$ with a 
simple pole along their respective graph;

\item consequently,  
$\hat s_\eps=\langle\mathbf x,\mathbf y\rangle^{-1}+\sum_c\langle\gt_c\mathbf x,
\mathbf y\rangle^{-1}+T$ is meromorphic on $\mathbb U$ with simple poles along the graphs 
$\mathbb P_1$, $\mathbb P_c$, and  
 holomorphic on $\mathbb U\setminus\bigcup\mathbb P$.

\end{enumerate}

\item {\rm Global convergence, invariance, meromorphicity.}

Let us assume (H) and $g\ge 2$.
 For $\mathbf x$, $\mathbf y\in\C^2\setminus0$ over $\Omega_0$ with 
$[\mathbf y]\notin\Gamma[\mathbf x]$ the series \eqref{eqshat} converges absolutely, and

\begin{equation}
\label{eqinvar}
\hat s_\eps(\tilde\lambda\mathbf x,\tilde\mu\mathbf y)
=\hat s_\eps(\mathbf x,\mathbf y), \qquad \tilde\lambda, \; \tilde\mu\in\Gt. 
\end{equation}

The function $\hat s_\eps$ is meromorphic on the cone over $\Omega_0\times\Omega_0$ 
(the points $\infty$ and $\Gamma\infty$ included), with polar set 
$\{[\mathbf y]\in\Gamma[\mathbf x]\}$, all poles simple. Equivalently, for fixed $y$ 
away from the orbit of the pole, $x\mapsto s_\eps(x,y)$ is meromorphic on $\Omega_0$ 
with simple poles at the points of $\Gamma y$. These project to the single point  
$[y]\in\Sg$, i.e., to the diagonal of $\Sg\times\Sg$. The only singularities in any chart 
are these poles: the point $x=\infty$ is not singular for $y\notin\Gamma\infty$, and there 
are no coordinate-dependent apparent singularities. For finite points $x$, $y$ with  
$\lambda x$, $\mu y\ne\infty$ the identity \eqref{eqinvar} reads
\begin{equation}\label{eqglue}
s_\eps(\lambda x,\mu y)\,j(\tilde\lambda,x)\,j(\tilde\mu,y)=s_\eps(x,y).
\end{equation}

\item {\rm Behaviour at $\infty$.} 
Set 
$s^\infty(y):=\hat s_\eps\bigl((1,0),(y,1)\bigr)
=1+\sum_{\gt\ne 1} ({c_{\gt}(\gamma\infty-y)})^{-1}$. 
It is meromorphic in $y\in\Omega_0\setminus\{\infty\}$, its poles are the points 
$\gamma\infty$, $\gamma\ne 1$, all simple, with residue $-1/c_{\gt}$ at $\gamma\infty$. 
For $y\in\Omega_0\setminus\Gamma\infty$,
\[
s_\eps(x,y)=\frac{s^\infty(y)}x+O\bigl(|x|^{-2}\bigr), \qquad |x|\to\infty, 
\]
locally uniformly in $y\in\Omega_0\setminus\Gamma\infty$. In particular 
$s_\eps(x,y)=O(1/|x|)$ 
 uniformly for $y$ in compact subsets of $\Omega_0\setminus\Gamma\infty$, 
and $S_\eps(\cdot,y)$ is regular at $x=\infty$ in the coordinate $u=1/x$ for such $y$. 

\end{enumerate}
\end{theorem}
 
Note that in $(i)$ no holomorphy across the polar graphs is claimed. For fixed $p$ 
the statement  
is meromorphicity on the interior of $\widehat{\mathcal K}_p\times\widehat{\mathcal K}_p$, 
 depending holomorphically on $p$, locally uniformly on $\mathfrak P$.

For $y\in\Gamma\infty$, $(iii)$ fails. If $y=\gamma\infty$, $\gamma\ne 1$, 
then $\lim_{|x|\to\infty}s_\eps(x,y)=-c_{\gt}\ne 0$. For instance,  
for $y=w_{-c}=\gamma_c(\infty)$ the length-one summand equals 
$j(\gt_c,x)/(\gamma_cx-w_{-c})=\lambda_c/\rho_c=-1/\lambda_c$, independent of $x$. 
Neither is the bound $O(1/|x|)$ uniform as $y\to\gamma\infty$, 
 since $x\,s_\eps(x,y)\to s^\infty(y)$ and $|s^\infty(y)|\to\infty$. 
The pairs $(\infty,\gamma\infty)$ are points of the polar set $\{y\in\Gamma x\}$, 
i.e., lie over the diagonal of $\Sg\times\Sg$.

\begin{proof}

(i) Lemma \ref{lemchart} bounds each of the $2g(2g-1)^{n-1}$ words of length $n\ge 2$ by 
$C_*\varpi^{n-2}$. Since $q:=(2g-1)\varpi\le\frac14$ the sum is at most 
$2g(2g-1)C_*\sum_{m\ge 0}q^m\le\frac83g(2g-1)C_*$, and $C_*,\varpi$ are continuous on 
the compact set $U$. Every term is holomorphic in $(p,\mathbf x,\mathbf y)$ 
(the entries of $\gt$ are polynomial in $w_c,\lambda_c^{\pm1}$) and, by Lemma \ref{lemchart}, 
nowhere zero in the denominator on $\widehat{\mathcal K}_p\times\widehat{\mathcal K}_p$.  
 The Weierstrass theorem gives holomorphy of the sum. The statements on the terms 
$n=0$, $1$ and on the polar set are immediate: 
$\langle\mathbf x,\mathbf y\rangle=0$ iff $\beta=\alpha$, and 
$\langle\gt_c\mathbf x,\mathbf y\rangle=0$ iff $\beta=\gamma_c\alpha$.

(ii) Uniformity from (i) is used only on the compact set 
$\widehat{\mathcal K}_p\times\widehat{\mathcal K}_p$ for fixed $p$ (or $p\in U$ compact).  
 The reindexing below is applied only to series that converge absolutely there, and no 
limiting operation is exchanged with the summation otherwise. 
Let $\alpha$, $\beta\in\Omega_0$, $\beta\notin\Gamma\alpha$. 
By the standard facts of Section \ref{secspin} ($\Omega_0=\Gamma\overline{\mathcal D}$) 
write $\alpha=\lambda\alpha_0$, $\beta=\mu\beta_0$ with $\alpha_0$, 
$\beta_0\in\overline{\mathcal D}\subset\operatorname{int}\widehat{\mathcal K}$ and 
$\lambda$, $\mu\in\Gamma$. 
Let $\tilde\lambda$, $\tilde\mu\in\Gt$ be the lifts and $\mathbf x_0$, $\mathbf y_0$ 
normalized representatives of $\alpha_0$, $\beta_0$. 
 Then $\beta_0\notin\Gamma\alpha_0$ 
(otherwise $\beta=\mu\beta_0\in\mu\Gamma\alpha_0=\Gamma\alpha$). For every $\gt\in\Gt$, 
by $\mathrm{SL}_2(\C)$-invariance of $\langle\cdot,\cdot\rangle$,
\[
\langle\gt\tilde\lambda\mathbf x_0,\tilde\mu\mathbf y_0\rangle
=\langle\tilde\mu^{-1}\gt\tilde\lambda\,\mathbf x_0,\mathbf y_0\rangle,
\]
and $\gt\mapsto\tilde\kappa:=\tilde\mu^{-1}\gt\tilde\lambda$ is a bijection of $\Gt$ onto 
itself (left and right translations). Thus, the series for 
$\hat s_\eps(\tilde\lambda\mathbf x_0,\tilde\mu\mathbf y_0)$ is a rearrangement of the series 
for $\hat s_\eps(\mathbf x_0,\mathbf y_0)$, which converges absolutely by (i)
 (together with the finitely many terms $n\le 1$, all finite as $\beta_0\notin\Gamma\alpha_0$). 
 Hence, we have absolute convergence and 
$\hat s_\eps(\tilde\lambda\mathbf x_0,\tilde\mu\mathbf y_0)
=\hat s_\eps(\mathbf x_0,\mathbf y_0)$. As every pair of vectors over 
$\Omega_0$ is $(t\tilde\lambda_1\mathbf x_0,t'\tilde\mu_1\mathbf y_0)$ and 
$\hat s_\eps$ is bihomogeneous of degree $(-1,-1)$, this proves \eqref{eqinvar} 
(i.e., replace $\tilde\lambda_1$, $\tilde\mu_1$ by 
$\tilde\lambda\tilde\lambda_1$, $\tilde\mu\tilde\mu_1$).

 \emph{Meromorphicity.} Near $(\alpha,\beta)$ write 
$(\mathbf x,\mathbf y)=(\tilde\lambda\mathbf x',\tilde\mu\mathbf y')$ 
with $(\mathbf x',\mathbf y')$ in a neighbourhood of $(\mathbf x_0,\mathbf y_0)$ 
contained in 
$\operatorname{int}\widehat{\mathcal K}\times\operatorname{int}\widehat{\mathcal K}$. 
 By \eqref{eqinvar} and (i), 
$\hat s_\eps(\mathbf x,\mathbf y)=\hat s_\eps(\mathbf x',\mathbf y')$ is meromorphic there. 

\emph{Poles.} The fixed points of the (loxodromic) elements $\gamma\ne 1$ belong to the limit 
set, 
thus inside $\Omega_0$ the graphs $\{y=\gamma x\}$, $\gamma\in\Gamma$, are pairwise disjoint. 
 By local finiteness of the family of graphs 
 (i.e., proper discontinuity, Section \ref{secspin}: 
only finitely many $\gamma$ have $\gamma K\cap K'\ne\emptyset$ for given compacta  
$K$, $K'\subset\Omega_0$),
 any point of $\Omega_0\times\Omega_0$ has a neighbourhood meeting at 
most one of them. At a point of the graph $\{y=\gamma x\}$ one term of the series 
(that of $\gt$) is singular and has a simple pole. The remaining terms sum to a 
function holomorphic there (by the same rearrangement and (i)). 
 Thus, there is no cancellation 
and the pole is simple. For finite points, 
 the relations $\hat s_\eps((x,1),(y,1))=s_\eps(x,y)$ 
and $\tilde g(x,1)=j(\tilde g,x)^{-1}(gx,1)$ turn \eqref{eqinvar} into \eqref{eqglue} 
by bihomogeneity.

(iii) For $\gt\ne 1$ one has 
$\langle\gt(1,0),(y,1)\rangle=a_{\gt}-c_{\gt}y=c_{\gt}(\gamma\infty-y)$, 
which gives the series for $s^\infty$. By (ii) (i.e., the pair $([(1,0)],[(y,1)])$  
is polar iff $y\in\Gamma\infty$) it is meromorphic and has simple poles with residue  
$-1/c_{\gt}$ at $\gamma\infty$. 
 Let $y_0\in\Omega_0\setminus\Gamma\infty$. The point $([(1,0)],[(y_0,1)])$ is not polar.   
Thus, by (ii) the function $(u,y)\mapsto\hat s_\eps((1,u),(y,1))$ is holomorphic near 
$(0,y_0)$. As $(x,1)=x\,(1,u)$ with $u=1/x$, bihomogeneity gives  
 $s_\eps(x,y)=u\,\hat s_\eps((1,u),(y,1))=u\bigl(s^\infty(y)+O(u)\bigr)$, locally uniformly 
in $y$.
  Now let $y_0=\gamma\infty$, $\gt=\begin{pmatrix}a&b\\c&d\end{pmatrix}\ne 1$. Thus,  
$y_0=a/c$. Near $(0,y_0)$ the only singular summand of $\hat s_\eps((1,u),(y,1))$ is the one 
of $\gt$. The others sum to a function $h(u,y)$ holomorphic there. Indeed, 
 by the rearrangement of (ii) with $\lambda=1$, $\mu=\gamma$, in the variables 
$\mathbf x=(1,u)$ and $\mathbf y'=\gt^{-1}(y,1)$, which lie near 
$\bigl((1,0),(1,0)\bigr)\in\operatorname{int}\widehat{\mathcal K}
\times\operatorname{int}\widehat{\mathcal K}$, the summand of $\gt$ becomes 
the term $1/\langle\mathbf x,\mathbf y'\rangle$ of length $0$. 
 The terms of length $1$ are regular there 
(since $\gamma_c\infty=w_{-c}\ne\infty$), and the tail is holomorphic by (i). 
 A computation, using $a-cy_0=0$ and $ad-bc=1$, 
 gives $\langle\gt(1,u),(y_0,1)\rangle=(a-cy_0)+(b-dy_0)u=-u/c$, 
 so that $s_\eps(x,y_0)=u\bigl(-c/u+h(u)\bigr)=-c+uh(u)\to-c$ 
as $u\to0$. Finally, $x\,s_\eps(x,y)\to s^\infty(y)$ 
as $|x|\to\infty$ for $y\notin\Gamma\infty$, and $s^\infty$ has a pole at $\gamma\infty$.
\end{proof}

Next, we deal with equivariance and global descent. 

\begin{proposition}
\label{propequi}
For $\tilde\lambda$, $\tilde\mu\in\Gt$ and finite $x$, $y\in\Omega_0$ with 
$y\notin\Gamma x$ and $\lambda x$, $\mu y\ne\infty$,
\[
s_\eps(\lambda x,y)\,j(\tilde\lambda,x)=s_\eps(x,y),
\qquad s_\eps(x,\mu y)\,j(\tilde\mu,y)=s_\eps(x,y). 
\]
 More generally, in homogeneous form, 
$\hat s_\eps(\tilde\lambda\mathbf x,\mathbf y)=\hat s_\eps(\mathbf x,\mathbf y)
=\hat s_\eps(\mathbf x,\tilde\mu\mathbf y)$ at all points of $\Omega_0$, 
including $\infty$ and $\Gamma\infty$, and there is no additive defect $\chi[\gamma](x,y)$. 
Consequently, $S_\eps$ descends to a meromorphic section of $\Lg\boxtimes\Lg$ on $\Sg\times\Sg$, holomorphic off the diagonal $\Delta$, with a simple pole along $\Delta$ and no other singularity.
\end{proposition}

\begin{proof}

The invariance is \eqref{eqinvar} with $\tilde\mu=1$ or $\tilde\lambda=1$, and the 
finite-point forms are \eqref{eqglue}. 

\emph{Descent:} as explained at the beginning of 
this section, sections of $\Lg\boxtimes\Lg$ over an open subset of $\Sg\times\Sg$ are 
the $\Gt\times\Gt$-invariant functions on the cone over its preimage in 
$\Omega_0\times\Omega_0$ that are bihomogeneous of degree $(-1,-1)$. 
Hence, $\hat s_\eps$ is a section of $\Lg\boxtimes\Lg$ over the complement of the image 
of its polar set. That image is the diagonal $\Delta\subset\Sg\times\Sg$: 
$y\in\Gamma x$ means that $x$, $y$ have the same class in $\Sg$. 
  $\Omega_0\times\Omega_0\to\Sg\times\Sg$ is a local biholomorphism since $\Gamma$ 
acts freely and properly discontinuously. 
 Meromorphicity across $\Delta$ and simplicity of the pole were proven in Theorem 
\ref{thmconv} (ii) at every point of $\Omega_0\times\Omega_0$, 
 in particular, at the points lying over $x=\infty$ or $y=\infty$, 
where the scalar trivialization by $(dz)^{1/2}$ is not available and the homogeneous 
form is used. Note that we do not claim regularity of $s_\eps$ at $x=\infty$,  
  The pairs $(\infty,\gamma\infty)$ are ordinary points of the polar divisor, 
and Theorem \ref{thmconv} (iii) describes how $s_\eps$ behaves there.
\end{proof}

Here below we discuss the skew-symmetry, local form, and projective connection. 

\begin{proposition}
\label{propskew}\label{propdiag}

\begin{enumerate}[(i)]
\item[]

\item $\hat s_\eps(\mathbf y,\mathbf x)=-\hat s_\eps(\mathbf x,\mathbf y)$ for all 
$\mathbf x$, $\mathbf y$ over $\Omega_0$ with $[\mathbf y]\notin\Gamma[\mathbf x]$. 
In particular, $s_\eps(y,x)=-s_\eps(x,y)$ for finite points.

\item {\rm Canonical form:} with $\pi(x,y)=\frac{(dx\,dy)^{1/2}}{x-y}$  
in any local coordinate near $\Delta$, 
 $S_\eps-\pi$ is a holomorphic section of $(\Lg\boxtimes\Lg)(\Delta)$ 
near $\Delta$, and its restriction to $\Delta$ 
 is the global constant $0\in\Gamma(\Sg,\mathcal O_\Sg)=\C$. 

\smallskip
\item[]
{\rm Coordinate form:} in any local coordinate this reads 
$s_\eps(x,y)=\frac1{x-y}+c_1(y)(x-y)+O((x-y)^2)$. 
 The coefficient $c_0(y)$ of $(x-y)^0$  
which represents the intrinsic restriction above in that coordinate 
vanishes identically. 

\item Let $z$ and $w$ be two local coordinates, $w=w(z)$, $w'=dw/dz$, and let $c_1^{(z)}$, 
$c_1^{(w)}$ denote the coefficient $c_1$ of part (ii) computed in the coordinates $z$ and $w$ 
respectively. Then $c_1^{(z)}(y)=w'(y)^2c_1^{(w)}(w(y))+\frac1{12}\{w;z\}(y)$. 
 Thus, $R_\eps:=12c_1$ is a projective connection: $R^{(z)}_\eps=w'^2R^{(w)}_\eps+\{w;z\}$.
\end{enumerate}
\end{proposition}

Note that the restriction of $(\Lg\boxtimes\Lg)(\Delta)$ to $\Delta$ in $(ii)$ is 
an element of 
$\Gamma(\Sg,(\Lg\boxtimes\Lg)(\Delta)|_\Delta)$, and 
$(\Lg\boxtimes\Lg)(\Delta)|_\Delta\cong\Lg^{\otimes2}
\otimes N_\Delta\cong K_\Sg\otimes T\Sg\cong\mathcal O_\Sg$ 
canonically since $\Lg^{\otimes2}\cong K_\Sg$ and $N_\Delta\cong T\Sg$ for the diagonal in 
$\Sg\times\Sg$. 
  In the coordinate form, vanishing of the coefficient $c_0(y)$ 
 is what ``no constant term'' means. 
 It is coordinate-independent because the restriction to $\Delta$ described above 
is defined without reference to a coordinate.

\begin{proof}

(i) For $\gt\in\Gt$, 
$\langle\gt\mathbf x,\mathbf y\rangle=-\langle\mathbf y,
\gt\mathbf x\rangle=-\langle\gt^{-1}\mathbf y,\mathbf x\rangle$ 
using antisymmetry of $\langle\cdot,\cdot\rangle$ and 
$\mathrm{SL}_2(\C)$-invariance.
  Summing over $\gt$ and reindexing $\gt\mapsto\gt^{-1}$, which is allowed by the absolute 
convergence of Theorem \ref{thmconv} (ii), gives 
$\hat s_\eps(\mathbf x,\mathbf y)=-\hat s_\eps(\mathbf y,\mathbf x)$.

(ii) Every point of $\Sg$ has a finite representative $x_0\in\Omega_0\setminus\{\infty\}$.  
 The class of $\infty$ is also represented by $w_{-c}=\gamma_c(\infty)$. 
 Thus, by Proposition \ref{propequi} it suffices to work near a diagonal point 
$(x_0,x_0)$ of finite points, in the coordinate $z=x$. By Theorem \ref{thmconv} (ii) 
the only singular summand of the series near $(x_0,x_0)$ is the one with $\gt=1$.  
 Note that for $\gamma\ne 1$ one has $x_0\ne\gamma x_0$, 
since the action is free, and the graphs 
 $\{y=\gamma x\}$ are pairwise disjoint.  
  Thus, $s_\eps-\frac1{x-y}$ is holomorphic there. 
 Write $s_\eps=\frac1{x-y}+c_0(y)+c_1(y)(x-y)+\dots$. Skew-symmetry at order $(x-y)^0$ 
gives $c_0=-c_0$, thus $c_0=0$.

(iii) Since $S_\eps$ is a well-defined section of $\Lg\boxtimes\Lg$, 
$s^{(z)}(x,y)=A\,s^{(w)}(w(x),w(y))$ with $A=(w'(x)w'(y))^{1/2}$ 
 with branch fixed by continuity along the diagonal. 
 Set $\xi:=x-y$, $a:=w''(y)/w'(y)$, $b:=w'''(y)/w'(y)$, 
and expand at fixed $y$: 
\[
\frac{w'(x)}{w'(y)}=1+a\xi+\tfrac b2\xi^2+O(\xi^3),\qquad \frac{w(x)-w(y)}{w'(y)\,\xi}
=1+\tfrac a2\xi+\tfrac b6\xi^2+O(\xi^3).
\]
Since $(1+e)^{-2}=1-2e+3e^2-\dots$, the second expansion gives 
$\bigl(\frac{w(x)-w(y)}{w'(y)\xi}\bigr)^{-2}=1-a\xi
+\bigl(\tfrac34a^2-\tfrac b3\bigr)\xi^2+O(\xi^3)$, and therefore
\begin{align*}
\frac{w'(x)w'(y)}{(w(x)-w(y))^2}&=\frac1{\xi^2}\Bigl(1+a\xi+\tfrac b2\xi^2\Bigr)
\Bigl(1-a\xi+\bigl(\tfrac34a^2-\tfrac b3\bigr)\xi^2\Bigr)+O(\xi)\\
&=\frac1{\xi^2}\Bigl(1+\bigl(\tfrac b6-\tfrac{a^2}4\bigr)\xi^2\Bigr)+O(\xi)
=\frac1{\xi^2}+\frac16\{w;z\}(y)+O(\xi),
\end{align*}
where $\{w;z\}:=\frac{w'''}{w'}-\frac32\bigl(\frac{w''}{w'}\bigr)^2=b-\frac32a^2$. 
The coefficient of $\xi^{-1}$ vanishes because the two terms $\pm a\xi$ cancel in the product. 
Taking the square root, with the branch $A/(w(x)-w(y))\sim+1/\xi$,
\[
\frac{A}{w(x)-w(y)}=\frac1\xi\Bigl(1+\tfrac1{12}\{w;z\}(y)\,\xi^2+O(\xi^3)\Bigr)
=\frac1\xi+\frac1{12}\{w;z\}(y)\,\xi+O(\xi^2),
\]
with no $\xi^0$ term, in agreement with part (ii). 
 Now insert 
$s^{(w)}(w(x),w(y))=\frac1{w(x)-w(y)}+c_1^{(w)}(w(y))\,(w(x)-w(y))+O(\xi^2)$ 
(using part (ii) in the coordinate $w$) into $s^{(z)}=A\,s^{(w)}$ and use 
$A\,(w(x)-w(y))=w'(y)^2\xi+O(\xi^2)$ to get 
\[
\frac1\xi+c_1^{(z)}(y)\,\xi+O(\xi^2)=\frac1\xi+\Bigl(\tfrac1{12}\{w;z\}(y)+w'(y)^2c_1^{(w)}(w(y))\Bigr)\xi+O(\xi^2).
\]
Comparing the coefficients of $\xi$ gives the required transformation law, 
and multiplying by $12$ gives $R^{(z)}_\eps=w'^2R^{(w)}_\eps+\{w;z\}$.
\end{proof}

This is an ordinary projective connection. In this paper we do not 
consider its counterparts such as super-projective connection 
 that would need a supergeometric structure, cf. \cite{Ma, P}. 

\subsection{Nonsingularity, uniqueness, identification}

Let us define the canonical Szeg\H{o} kernel. 

\begin{definition} 
For a spin structure $L$, a Szeg\H{o} kernel is a meromorphic section $S_L$ of 
$L\boxtimes L$ on $\Sg\times\Sg$, holomorphic off $\Delta$, with 
$S_L-\frac{(dx\,dy)^{1/2}}{x-y}$ holomorphic near $\Delta$ in every local coordinate. 
\end{definition}

 This is the intrinsic meaning of ``simple pole with unit principal part''.   
 The polar part is the canonical section $1$ of 
$(L\boxtimes L)(\Delta)|_\Delta\cong L^2\otimes T\cong\mathcal O$. 
 It is coordinate independent by Proposition \ref{propdiag} (ii). 
We now distinguish: (a) the series of Definition \ref{defS}; 
 (b) the descended section $S_\eps$ of Proposition \ref{propequi}; and 
(c) a canonical kernel $S_L$ as above. 

\begin{remark}
\label{remdiagdom} 

Now consider the domain and invariance of the diagonal expansion. 
The expansion $s_\eps(x,y)=\frac1{x-y}+c_1(y)(x-y)+O((x-y)^2)$ of Proposition 
\ref{propdiag} (ii) holds for $y$ in a compact subset of $\Omega_0$ and $|x-y|$ smaller 
than the distance from $y$ to the other points of $\Gamma y$ 
(a positive, locally bounded-below quantity since $\Gamma$ acts properly discontinuously 
on $\Omega_0$). The coefficients $c_0\equiv0$ and $c_1$ are coordinate-dependent.  
 Under a coordinate change $y\mapsto\phi(y)$ the term $1/(x-y)$ is replaced by 
$\phi'(x)^{1/2}\phi'(y)^{1/2}/(\phi(x)-\phi(y))$, 
 whose expansion is $\frac1{x-y}+\frac{1}{12}\{\phi,y\}(x-y)+\dots$ 
 where $\{\phi,y\}$ is the Schwarzian.   
 $c_1$ changes by $\phi'^2c_1\circ\phi$ added by the Schwarzian term. 
 This transformation law is what makes $R_\eps=12c_1$ a projective connection 
(normalized so that $R\mapsto\phi'^2R\circ\phi+\{\phi,y\}$). On $\Sg$ the quantity $c_1$ 
is thus invariantly defined only as a projective connection. 
\end{remark}

The theorem below concerns the nonsingularity, uniqueness, and identification. 
 Let us assume (H) with $g\ge 2$ or, more generally, the hypotheses of Theorem \ref{thmE}. 
 Recall Proposition \ref{propequi}. The statements of the theorem below concern only 
 the $\Ng$ lifted spin structures $\Lg$.  
   In particular, we do not consider the remaining spin structures  
(cf.  Remark \ref{remdict}), for some of which $H^0\ne 0$ can occur.   
 The vanishing $H^0(\Sg,\Lg)=0$ is obtained in the proof from the existence of $S_\eps$, 
 a consequence of the convergence hypotheses.
 
\begin{theorem}
\label{thmnonsing}

\begin{enumerate}[(i)]

\item $H^0(\Sg,\Lg)=0$ for every one of the $\Ng$ lifts, for all parameters satisfying (H). 

\item If $S'$ is any section of $\Lg\boxtimes\Lg$ with the properties of a Szeg\H{o} kernel,
 then $S'=S_\eps$. Thus, (a)=(b)=(c) mentioned above. 
\end{enumerate}
\end{theorem}

$(i)$ is a statement about that family on this domain only. We do not claim that all  
even spin structures are nonsingular. 

\begin{proof}
(i) Let $L=\Lg$ and suppose $k:=h^0(L)>0$. Choose $0\ne\phi\in H^0(L)$. Its zero divisor  
has degree $\deg L=g-1$, hence is finite. We choose $p\in\Sg$ outside it. The sections  
of $L(-p)$ are the sections of $L$ vanishing at $p$, i.e., the kernel of the evaluation map 
$H^0(L)\to L|_p$, which has rank $\le 1$, and $\phi\notin H^0(L(-p))$. 
 Hence, $h^0(L(-p))=k-1$. Riemann-Roch for $L(p)$, of degree $g$, gives 
$h^0(L(p))-h^1(L(p))=1$. Serre duality with $K_\Sg\otimes L^{-1}\cong L$ gives 
$h^1(L(p))=h^0(L(-p))=k-1$. Hence, $h^0(L(p))=k=h^0(L)$, and the injection 
$H^0(L)\hookrightarrow H^0(L(p))$ is onto. Every section of $L(p)$ is regular at $p$.  
 On the other hand, take $y\in\Omega_0$ with $[y]=p$. By Proposition \ref{propequi}  
the kernel $S_\eps(\cdot,y)$ is a global meromorphic section of $L$ on $\Sg$ whose only 
singularity is a simple pole at $p$ with nonzero principal part 
(Proposition \ref{propdiag} (ii)), i.e., a section of $L(p)$ not regular at $p$. 
This contradiction shows $k=0$. 
 Proposition \ref{propequi} and Theorem \ref{thmconv} (ii) provide 
the global existence of $S_\eps(\cdot,y)$, including across the points over $\infty$. 

(ii) The difference $S'(\cdot,y)-S_\eps(\cdot,y)$ is holomorphic on $\Sg$ since 
 there is no pole at $x=y$ by the principal-part condition. Thus, it vanishes by (i).
\end{proof}

\begin{remark}
\label{remcompare}

 Let us make a comparison with the Bers case.  
For $N\ge 2$ the factor $\prod_{j\in\mathcal L}\frac{y-A_j}{x-A_j}$ is not invariant under 
joint M\"obius action (the $A_j$ do not move), which blocks the reindexing used above and 
produces a cocycle $\chi_N[\gamma]$ with values in $\mathfrak P_{2N-2}$. At $N=\half$ there 
are $2N-1=0$ auxiliary points, and the joint invariance of $\pi$ makes the reindexing 
legitimate by absolute convergence. There is no space of global half-differentials on 
$\mathbb P^1$ that could play the role of the space of polynomials of degree $2N-2$.
\end{remark}

\begin{remark}
\label{remlambda0}
   
At $\lambda_a=0$ one has $\rho_a=0$. Thus, $\gamma_a(z)\equiv w_{-a}$ is a constant map, 
the matrix $\gt_a$ is undefined, and neither $\Gamma$ nor $\Sg=\Omega_0/\Gamma$ is defined 
as above. $\lambda_a=0$ corresponds to a noded limit and is not a point of the parameter 
space $\mathfrak P$. We then formulate all statements are for $\lambda_a\ne 0$. 
The limit $\lambda_a\to 0$ is quantified, with an error bound, by Theorem \ref{thmexp}. 
\end{remark}

\begin{remark} 
\label{remcont}

$S_\eps$ is a function of $\lambda_a$, not of $\rho_a=-\lambda_a^2$. 
The sign of $\lambda_a$ labels the lift. Continuing $\lambda_a\mapsto e^{i\phi}\lambda_a$ 
along $\phi\in[0,\pi]$ assuming (H) (a half-turn around the excluded point $\lambda_a=0$) 
exchanges the two lifts of $\gamma_a$. Hence, it permutes the $\Ng$ kernels 
$\{S_\eps\}$. Continuation of $S_\eps$ beyond the domains of Section \ref{secE} 
(for instance, across the theta-null divisor $\theta[\delta](0)=0$, where by Theorem 
\ref{thmnonsing} the series cannot converge in the sense of Theorem \ref{thmE}) is not 
treated here (see Section \ref{secopen}).
\end{remark}

\subsection{Convergence beyond (H). M\"obius covariance}\label{secE}

Condition (H) is a coordinate-dependent, small neighbourhood of the maximal
 degeneration. The following facts show that the results of Section 3 
do not depend on it. The statement on the M\"obius covariance is the following. 
For $h\in\mathrm{SL}_2(\C)$ set a lift of $h\Gamma h^{-1}$ as 
 $\Gt^h:=h\Gt h^{-1}$. Then 

\begin{proposition}
\label{propmob}

$\hat s_{\Gt^h}(h\mathbf x,h\mathbf y)=\hat s_{\Gt}(\mathbf x,\mathbf y)$ 
whenever one side converges absolutely.
\end{proposition}

 In particular, the $3g$ parameters $(w,\lambda)$ are redundant by the action 
of $\mathrm{PSL}_2(\C)$, and the kernels $S_\eps$ are objects attached to a point of
 Schottky space together with a lift.

\begin{proof}
$\gt\mapsto h\gt h^{-1}$ is a bijection $\Gt\to\Gt^h$, and 
$\langle h\gt h^{-1}h\mathbf x,h\mathbf y\rangle=\langle h\gt\mathbf x,h\mathbf y\rangle
=\langle\gt\mathbf x,\mathbf y\rangle$ by $\mathrm{SL}_2(\C)$-invariance of
 $\langle\cdot,\cdot\rangle$.
\end{proof}

Let $\Nrm{\cdot}$ be the operator norm on $\C^2$ (Euclidean). We consider the condition
\begin{equation}\tag{E}
\sum_{\gt\in\Gt}\Nrm{\gt}^{-1}<\infty .
\end{equation}
It does not depend on the choice of lift since the lifts differ by signs, and it is invariant
 under conjugation in $\mathrm{SL}_2(\C)$ (norms change by bounded factors), hence intrinsic 
to $\Gamma$.
Under conditions (D) and (E) the convergence is described by 

\begin{theorem}
\label{thmE}

Let $\Gamma$ be a Schottky group of genus $g\ge 1$ as in Section \ref{secspin}, 
satisfying (D), and let $\Gt$ be a lift satisfying (E). Then:

\begin{enumerate}[(i)]

\item For compact $K,K'\subset\Omega_0$ there are a finite set $F\subset\Gt$ and $\eta>0$
 such that $|\langle\gt\mathbf x,\mathbf y\rangle|\ge\eta\Nrm{\gt}$ for all $\gt\notin F$ 
and all unit vectors $\mathbf x$, $\mathbf y$ over $K$, $K'$. Hence, the series \eqref{eqshat} 
converges absolutely and locally uniformly on 
$\{(\mathbf x,\mathbf y)\ \text{over}\ \Omega_0\times\Omega_0:\
 [\mathbf y]\notin\Gamma[\mathbf x]\}$ (the points $\infty$ and $\Gamma\infty$ included), 
is meromorphic on the cone over $\Omega_0\times\Omega_0$, with simple poles along 
$\{[\mathbf y]\in\Gamma[\mathbf x]\}$, and satisfies the invariance \eqref{eqinvar} and the 
gluing relation \eqref{eqglue}.

\item Consequently, Proposition \ref{propequi}, Proposition \ref{propskew} (i)-(iii), 
Theorem \ref{thmnonsing} (in particular $H^0(\Sg,\Lg)=0$), Propositions \ref{propfay},
 \ref{propS2} and Corollary \ref{corberg}, Theorem \ref{thmred}, Proposition \ref{proppf},
 Proposition \ref{propT} hold, with the same proofs, for every such $\Gt$.

\item For every $g\ge 1$, condition (H) implies (D) and (E).  
 If the critical exponent $\delta(\Gamma)<\half$ 
 (the abscissa of convergence of $\sum_\gamma e^{-s\,d(o,\gamma o)}$ for the hyperbolic 
distance in $\mathbb H^3$), then (E) holds.
\end{enumerate}

\end{theorem}

Note that condition (H) is meaningful for all $g\ge 1$. In the proof 
the constants of Lemma \ref{lemchart} are not used. Instead, we use only Lemma \ref{lemword},
 valid for all $g\ge 1$. Genus one does not need (H): it is covered by (E)  
 (Proposition \ref{propg1}).

\begin{proof}

(i) Fix $K$, $K'$ and $\eta_0>0$ so small that the closed chordal $\eta_0$-neighbourhood 
$K''$ of $K'$ belongs to $\Omega_0$ ($\Omega_0$ is open). By proper discontinuity the set 
$F_1:=\{\gamma:\gamma K\cap K''\ne\emptyset\}$ is finite. For $\gamma\notin F_1$, $\gamma x$ 
has chordal distance $\ge\eta_0$ from $y$ for $x\in K,y\in K'$. For unit vectors one has 
$|\langle\mathbf u,\mathbf v\rangle|=|\mathbf u||\mathbf v|\,\chi([\mathbf u],[\mathbf v])$
 with $\chi\in[0,1]$ the chordal distance (Lagrange identity 
$|u_1v_2-u_2v_1|^2=|\mathbf u|^2|\mathbf v|^2-|u_1\bar v_1+u_2\bar v_2|^2$). Thus,  
$|\langle\gt\mathbf x,\mathbf y\rangle|\ge\eta_0\,|\gt\mathbf x|$ for $\gamma\notin F_1$.
Next, let $\sigma_1\ge\sigma_2=\sigma_1^{-1}$ be the singular values of $\gt$ ($\det\gt=1$), 
$\mathbf v_1,\mathbf v_2$ orthonormal right singular vectors. Denoting  
$\mathbf x=c_1\mathbf v_1+c_2\mathbf v_2$, we obtain 
$|\gt\mathbf x|^2=\sigma_1^2|c_1|^2+\sigma_2^2|c_2|^2\ge\Nrm{\gt}^2\chi([\mathbf x],
[\mathbf v_2])^2$ again by Lagrange identity, $|c_1|=\chi([\mathbf x],[\mathbf v_2])$. 

 Consider the followin claim: there are a finite $F_2$ and $\eta_1>0$ with 
$\chi([\mathbf x],[\mathbf v_2(\gt)])\ge\eta_1$ for $\gt\notin F_2$, $\mathbf x$ over $K$. 
Otherwise there are distinct $\gt_n$ and $x_n\in K$ with 
$\chi(x_n,[\mathbf v_2(\gt_n)])\to0$. Passing to a subsequence, $x_n\to b\in K$, 
$[\mathbf v_2(\gt_n)]\to b$ and $\gt_n/\Nrm{\gt_n}\to A=\mathbf u_1\mathbf v_1^*$ of rank one 
(here $\Nrm{\gt_n}\to\infty$ since $\Gt$ is discrete), with $\ker A=b$. For $y\ne b$ we get 
$\gamma_ny\to a:=[\mathbf u_1]$, with the $\gamma_ny$ pairwise distinct points of an orbit of
 a point of $\Omega_0$, hence $a\in\Lambda$. Applying the same reasoning to $\gt_n^{-1}$ 
(whose dominant left singular vector is $\mathbf v_2(\gt_n)$) gives $b\in\Lambda$. 
This contradicts $b\in K\subset\Omega_0$. Set $F:=F_1\cup F_2$ and $\eta:=\eta_0\eta_1$. 
 Then $|\langle\gt\mathbf x,\mathbf y\rangle|\ge\eta\Nrm{\gt}$ for $\gt\notin F$, and (E) 
gives uniform absolute convergence of the sum over $\gt\notin F$ on $K\times K'$. 
 The finitely many remaining terms are meromorphic. Bihomogeneity gives the same for vectors
  of bounded norm. Invariance \eqref{eqinvar} follows by reindexing 
$\gt\mapsto\tilde\mu^{-1}\gt\tilde\lambda$, which is allowed by absolute convergence,
 and \eqref{eqglue} by the computation in Theorem \ref{thmconv} (ii). 

\emph{Poles:} at a point $(\alpha,\gamma_0\alpha)$ only the term of $\gt_0$ is singular 
(due to freeness of the action). It has a simple pole, and the sum of the others 
is holomorphic nearby (finitely many terms of $F$ with $\gamma\ne\gamma_0$, 
and others are the uniformly convergent). 

(ii) The proofs of the listed statements use only the invariance \eqref{eqinvar}, 
the meromorphicity of $\hat s_\eps$ on the cone over $\Omega_0\times\Omega_0$ with simple 
poles  along $\{y\in\Gamma x\}$ and local finiteness of the family of graphs. All of these 
are now provided by (i). For $g=1$ the Riemann-Roch argument of Theorem \ref{thmnonsing} 
applies: $\deg\Lg=0$ and the zero divisor of a nonzero section is empty. 

(iii) That (H) implies (D) was noted above. Let $x_0:=w_1+d/2$. Then 
$|x_0-w_c|\ge d/2\ge r_c/2$ for all $c$, i.e., $x_0\in\mathcal K$. By Lemma \ref{lemword}, 
$|c_{\gt}x_0+d_{\gt}|=|j(\gt,x_0)|^{-1}\ge\frac12\varpi^{-(n-1)}$ for $|\gt|=n\ge 1$, and 
$\Nrm{\gt}\ge|\gt(x_0,1)|/|(x_0,1)|\ge|c_{\gt}x_0+d_{\gt}|/(1+|x_0|)$. Therefore,  
$\sum_{\gt\ne 1}\Nrm{\gt}^{-1}\le 2(1+|x_0|)
\sum_{n\ge 1}2g(2g-1)^{n-1}\varpi^{n-1}\le 4g(1+|x_0|)\sum_{m\ge 0}4^{-m}<\infty$. 
For the last claim, $\Nrm{\gt}^2\asymp e^{d(o,\gamma o)}$ for a base point $o\in\mathbb H^3$, 
thus $\sum_\gamma\Nrm{\gt}^{-1}\asymp\sum_\gamma e^{-d(o,\gamma o)/2}$, which converges when 
$\tfrac12>\delta(\Gamma)$ (definition of the critical exponent, see \cite{Su}).
\end{proof}

A few remarks are in order. 
(1) Numerically (genus two, random data with $r/d$ up to $0.3$, far outside (H)),  
the series converges geometrically in the word length and the relation \eqref{eqglue} 
holds to the accuracy of the truncation. 
 (2) By Theorem \ref{thmnonsing} and the classical criterion $\theta[\delta](0)=0$ iff $h^0>0$, 
(E) fails at every Schottky group for which $\theta[\delta(\eps)](0)=0$. 
Thus the locus where (E) holds avoids those theta-null divisors. (3) Whether (E) 
 fails when $\delta(\Gamma)\ge\half$ is not claimed.

\subsection{Relation to the prime form and the canonical bidifferential}

\textit{Normalization:} fix the marking $\alpha_a=\partial D_a$, $\beta_a$ a path joining 
$\partial D_{-a}$ to $\partial D_a$ through $\mathcal D$. Let $\zeta=(\zeta_1,\dots,\zeta_g)$ 
be the holomorphic $1$-forms normalized by $\oint_{\alpha_a}\zeta_b=\delta_{ab}$, 
$\tau_{ab}=\oint_{\beta_a}\zeta_b$, thus $\mathrm{Im}\,\tau>0$, and
\[
\theta[\delta](z)=\sum_{n\in\mathbb Z^g}\exp\bigl(\pi i(n+\delta')^T\tau(n+\delta')+2\pi i(n+\delta')^T(z+\delta'')\bigr).
\]
Let $B(x,y)$ be the canonical bidifferential normalized by $\oint_{\alpha_a}B(\cdot,y)=0$. 
 We denote by $E$ the prime form, $E(x,y)=(x-y)(dx)^{-1/2}(dy)^{-1/2}(1+O((x-y)^2))$.  
   Thus, $B(x,y)=\dd_x\dd_y\log E\,dx\,dy=\frac{dx\,dy}{(x-y)^2}+\dots$ 
(note that $\dd_x\dd_y\log(x-y)=+(x-y)^{-2}$). Since $E$ is skew symmetric, 
$E(x,y)=(x-y)(dx\,dy)^{-1/2}\bigl(1+e_2(y)(x-y)^2+O((x-y)^3)\bigr)$, 
i.e., there is no linear term.   
 We use the following standard facts \cite{Fa, Mu}. Note that $B=\dd_x\dd_y\log E$. 
 As a function of $x$ in a fixed local coordinate, $E(x,y)$ is single valued along 
$\alpha$-cycles and  
is multiplied by $-\exp(-\pi i\tau_{bb}-2\pi i\int_y^x\zeta_b)$ along $\beta_b$.  
 For $\theta[\delta](0)\ne 0$, the Szeg\H{o} kernel of $L_\delta$ is the expression 
in Proposition \ref{propfay}.  
 The following proposition is due to Fay \cite{Fa, Mu}. 

\begin{proposition}
\label{propfay}

Assume the hypotheses of Theorem \ref{thmE} (in particular (H) with $g\ge 2$, or $g=1$ and (E)).

\begin{enumerate}[(a)]

\item For $\delta=\delta(\eps)$ (Convention of Section 2) one has $\theta[\delta](0)\ne 0$ 
(by Theorem \ref{thmnonsing} (i) and the criterion above), and
\[
S_\eps(x,y)=\frac{\theta[\delta]\!\bigl(\int_y^x\zeta\bigr)}{\theta[\delta](0)\,E(x,y)}.
\]

\item $\delta'=0$, i.e., $\delta(\eps)=[0;\delta'']$ with $\delta''\in\{0,\half\}^g$. 
In particular $\delta(\eps)$ is even. 

\end{enumerate}
For odd characteristics (not in our family) $\theta[\delta](0)=0$ and the formula 
is unavailable.

\end{proposition}

\begin{proof}

(a) By the quoted fact the right-hand side is a Szeg\H{o} kernel of $L_\delta\cong\Lg$, 
and by Theorem \ref{thmnonsing} (ii) it equals $S_\eps$. (b) Let $x$ run once around 
$\alpha_b=\partial D_b$ in $\Omega_0$ with $y$ fixed. The scalar function $s_\eps(\cdot,y)$
 is single valued on $\Omega_0$ (the loop is closed in $\Omega_0$, where the trivialization by 
$(dz)^{1/2}$ is global). On the right-hand side, $\int_y^x\zeta_c$ increases by $\delta_{bc}$,
 thus $\theta[\delta]$ is multiplied by $e^{2\pi i\delta'_b}$, while $E(x,y)$ is single valued 
along $\alpha_b$. Hence $e^{2\pi i\delta'_b}=1$ for all $b$, i.e., $\delta'=0$. 
 Then $4\delta'\cdot\delta''=0$ and $\delta$ is even.
\end{proof}
 
 The identification of the bundle with $L_{\delta(\eps)}$ is the Convention of Section 2.
   Let us check normalization. 
As $y\to x$, $\theta[\delta]\bigl(\int_y^x\zeta\bigr)\to\theta[\delta](0)$ 
and $E(x,y)\to(x-y)(dx)^{-1/2}(dy)^{-1/2}$ by the stated normalization of $E$. 
 Thus, the right-hand side tends to 
\[
\dfrac{\theta[\delta](0)}{\theta[\delta](0)(x-y)(dx)^{-1/2}(dy)^{-1/2}}
=\dfrac{(dx\,dy)^{1/2}}{x-y},
\]
 which is the principal part of $S_\eps$
fixed in Definition \ref{defS} with no extra sign,
 confirming that the stated normalization of $E$ is the one compatible with our $S_\eps$.

\subsection{Review of conventions}

 (a) \emph{Characteristics:} with the theta function displayed above, 
$\theta[\delta](z+e_b)=e^{2\pi i\delta'_b}\theta[\delta](z)$ for the $b$-th unit vector 
$e_b$ (periods along $\alpha_b$). This is the multiplier used in Proposition 
\ref{propfay} (b).  
  Only the intrinsic bundle isomorphism $L_{\delta(\eps)}\cong\Lg$ is used 
for $\delta(\eps)$ with no sign convention for the lifts entering.

 (b) \emph{Parity and sign:} for even $\delta$, $\theta[\delta](-z)=\theta[\delta](z)$ 
and $E(y,x)=-E(x,y)$. Thus, the right-hand side of Proposition \ref{propfay} (a) 
is skew symmetric in $(x,y)$ as required by Proposition \ref{propdiag} (i). 

(c) \emph{Overall sign:} the sign of the prime form is fixed here by 
$E(x,y)\sim(x-y)(dx)^{-1/2}(dy)^{-1/2}$, equivalently by the requirement that 
$S_\eps$ have principal part $+\frac{(dx\,dy)^{1/2}}{x-y}$ (Definition \ref{defS}). 
The footnote in \cite{TZ1}, Section 2, mentiones that their prime form differs by a 
factor $-1$ from that of \cite{Mu, Fa}. With the opposite sign for $E$ one obtains 
the same kernel with the opposite sign in the formula.  

(d) The formula for $S_\eps$ is also confirmed independently at genus one 
(Proposition \ref{propg1}), including the sign and the normalization of the prime form.

\begin{proposition}
\label{propS2}

$S_\eps(x,y)^2-B(x,y)=\sum_{i,j}c_{ij}\zeta_i(x)\zeta_j(y)$
 is a holomorphic bidifferential. 

\end{proposition}

\begin{proof}

The stetement comes from 
  using $H^0(\Sg\times\Sg,K\boxtimes K)\cong H^0(\Sg,K)\otimes H^0(\Sg,K)$, 
a standard consequence of the K\"unneth formula, 
together with the choice of basis $\zeta_i$. 
 $S_\eps^2$ is a section of $K\boxtimes K$ (signs square to $1$ 
by Proposition \ref{propequi}) with principal part $\frac{dx\,dy}{(x-y)^2}$ near 
$\Delta$ and, by Proposition \ref{propdiag} (ii), there is no term of order $(x-y)^{-1}$. 
The same holds for $B$ by definition of the canonical bidifferential. Hence, $S_\eps^2-B$ 
is holomorphic on $\Sg\times\Sg$, and the displayed decomposition follows from the K\"unneth 
 isomorphism. 
\end{proof}

Let us describe now projective connection and Bergman projective connection. 

\begin{corollary}
\label{corberg}
Let $S_B$ be the Bergman projective connection, defined by 
\[
B(x,y)=\frac{dx\,dy}{(x-y)^2}+\frac16S_B(y)\,dx\,dy+O(x-y) 
\]
in a local coordinate. Then in every local coordinate
\[
R_\eps=12\,c_1=S_B+6\sum_{i,j}c_{ij}\,\zeta_i\zeta_j,
\]
with $c_{ij}$ as in Proposition \ref{propS2}. 
 Moreover $c_{ij}=\dd_i\dd_j\log\theta[\delta](0)$, $\delta=\delta(\eps)$, where 
$\dd_i=\partial/\partial z_i$ and $z_i=\int\zeta_i$, at least as quadratic differentials, 
i.e., $\sum c_{ij}\zeta_i\zeta_j=\sum\dd_i\dd_j\log\theta[\delta](0)\,\zeta_i\zeta_j$.

\end{corollary}

\begin{proof}

By Proposition \ref{propdiag} (ii), $S_\eps^2=(x-y)^{-2}+2c_1(y)+O(x-y)$ 
(the cross term $2c_0/(x-y)$ vanishes). 
 $B$ has the expansion above. Subtracting, we obtain 
 $\sum c_{ij}\zeta_i(y)\zeta_j(y)=2c_1(y)-\frac16S_B(y)$, 
which is the first formula. 
 For the second, set $\xi=x-y$ and use Proposition \ref{propfay} 
(a). Since $\delta$ is even, 
$\theta[\delta](\int_y^x\zeta)/\theta[\delta](0)
=1+\frac12\sum\dd_i\dd_j\log\theta[\delta](0)\,\zeta_i(y)\zeta_j(y)\,\xi^2+O(\xi^3)$, 
 and $E^{-2}=\xi^{-2}(dx\,dy)(1-2e_2\xi^2+\dots)$ while 
$B=\dd_x\dd_y\log E=\xi^{-2}-2e_2+O(\xi)$, thus $E^{-2}-B=O(\xi)$. 
Hence, $S_\eps^2-B=E^{-2}\bigl((\theta[\delta]/\theta[\delta](0))^2-1\bigr)
+(E^{-2}-B)=\sum\dd_i\dd_j\log\theta[\delta](0)\zeta_i(y)\zeta_j(y)+O(\xi)$.
\end{proof}
  Since $S_B$ is a projective connection, this is consistent with Proposition 
\ref{propdiag} (iii). For $g=1$ it is checked in Proposition \ref{propg1}(c) below.

\subsection{Genus one}

Let $g=1$ and let $\Gt=\langle\gt\rangle$ with $\gt$ loxodromic of eigenvalues $\kappa^{\pm1}$,
 $0<|\kappa|<1$. This is the general form of a lifted genus one Schottky group.  
(D) is empty and (E) holds because $\Nrm{\gt^n}\asymp|\kappa|^{-|n|}$.  
 Set $q=\kappa^2$ and fix $\tau$ with $e^{2\pi i\tau}=q$, $\mathrm{Im}\,\tau>0$. 
By Proposition \ref{propmob} we may conjugate so that 
$\gt=\mathrm{diag}(\kappa,\kappa^{-1})$, i.e., $\gamma x=qx$, $\Omega_0=\C^\times$, 
$\Sg=\C^\times/q^{\Z}$. Then we have 

\begin{proposition}
\label{propg1}

\begin{enumerate}[(a)]

\item For $x$, $y\in\C^\times$, $y\notin q^{\Z}x$:  
$\;s_\eps(x,y)=\sum_{n\in\Z}\dfrac{\kappa^n}{q^nx-y}$.

\item Define $\delta''\in\{0,\half\}$ by 
$\kappa=-e^{2\pi i\delta''}e^{i\pi\tau}$ (note that $(\kappa e^{-i\pi\tau})^2=1$). 
Then with $x=e^{2\pi iu}$, $y=e^{2\pi iv}$,
\[
s_\eps(x,y)=\frac1{2\pi i\,e^{\pi i(u+v)}}
\cdot\frac{\theta[0;\delta''](u-v)\,\theta_1'(0)}{\theta[0;\delta''](0)\,
\theta_1(u-v)},\qquad\theta_1:=\theta[\tfrac12;\tfrac12],\ \ '=\tfrac d{du}, 
\]
where the right side is independent of the branches of $u$, $v$.

\item In the coordinate $x$, 
$R_\eps(y)=y^{-2}\Bigl[\tfrac12-\tfrac1{4\pi^2}
\Bigl(6\,\tfrac{\theta[0;\delta'']''(0)}
{\theta[0;\delta''](0)}-2\,\tfrac{\theta_1'''(0)}{\theta_1'(0)}\Bigr)\Bigr]$, 
and $S_B=-2\theta_1'''(0)/\theta_1'(0)$ in the coordinate $u$.
\end{enumerate}
\end{proposition}

Note that it follows from (b) that    
 Proposition \ref{propfay} holds with $\delta(\eps)=[0;\delta'']$, 
and $\kappa\mapsto-\kappa$ (the other lift) exchanges $\delta''=0$ 
and $\delta''=\half$. Replacing $\tau$ by $\tau+1$ 
(another marking, $\beta\mapsto\beta+\alpha$) exchanges the two values as well. 

\begin{proof}
(a) $\langle\gt^n(x,1),(y,1)\rangle=\kappa^nx-\kappa^{-n}y$. 

(b) Let $\Phi(x,y)$ be the right-hand side.
  Using $\theta[0;\delta''](z+\tau)=e^{-\pi i\tau-2\pi i(z+\delta'')}\theta[0;\delta''](z)$, 
$\theta_1(z+\tau)=-e^{-\pi i\tau-2\pi iz}\theta_1(z)$, the one-periodicity of 
$\theta[0;\delta'']$ and the sign change of $\theta_1$ under $z\mapsto z+1$ 
(compensated by $e^{\pi i(u+v)}\mapsto-e^{\pi i(u+v)}$), one finds that $\Phi$ is 
well defined on $\C^\times\times\C^\times$. Also, 
$\Phi(qx,y)=-e^{-2\pi i\delta''}e^{-i\pi\tau}\Phi(x,y)=\kappa^{-1}\Phi(x,y)$ and 
likewise $\Phi(x,qy)=\kappa^{-1}\Phi(x,y)$, which is the gluing relation \eqref{eqglue} 
($j(\gt,x)=\kappa$). $\Phi$ has simple poles on $y\in q^{\Z}x$, and near $y=x$, 
$\theta_1(u-v)/\theta_1'(0)=(x-y)/(2\pi iy)(1+O(x-y))$. Thus, $\Phi=\frac1{x-y}+O(1)$. 
Now fix $y$. Then $D(x):=s_\eps(x,y)-\Phi(x,y)$ has no pole at $x=y$ and, by the common 
automorphy law, none at $q^nx$. Thus, $D$ is holomorphic on $\C^\times$ with 
$D(qx)=\kappa^{-1}D(x)$. Writing $D=\sum a_mx^m$ gives $a_m(\kappa^{2m+1}-1)=0$, 
and $\kappa^{2m+1}\ne 1$ for every integer $m$ ($0<|\kappa|<1$), thus $D=0$.  

(c) From (b), in the coordinate $u$: 
$c_1^{(u)}
=\frac12\theta[0;\delta'']''(0)/\theta[0;\delta''](0)-\theta_1'''(0)/(6\theta_1'(0))$ 
and $R^{(u)}=12c_1^{(u)}$. 
 $B=-\partial_u^2\log\theta_1(u-v)$ gives $S_B^{(u)}=-2\theta_1'''(0)/\theta_1'(0)$. 
For $u=\log x/(2\pi i)$, $u'=1/(2\pi ix)$ and $\{u;x\}=\frac1{2x^2}$, 
thus Proposition \ref{propdiag} (iii) gives 
$R^{(x)}=-\frac{1}{4\pi^2x^2}R^{(u)}+\frac1{2x^2}$, which is the required formula.
\end{proof}
 
 In particular, the genus one kernel for $\delta''=0$ and $\delta''=\half$ 
is the classical torus Szeg\H{o} kernel of the two even characteristics with trivial 
$\alpha$-holonomy. The other two characteristics of the torus are not in the family.

\begin{remark}
\label{remdict}

 For $g\ge 2$ we have the following. 
Let $\kappa_a$ be the eigenvalue of $\gt_a$ of modulus $<1$. 
Then $\kappa_a=-\lambda_a/(w_a-w_{-a})\,(1+O(r/d))$. At genus one, 
(b) gives $\kappa=-e^{2\pi i\delta''}e^{i\pi\tau}$. 
 For $g\ge 2$ we expect, for the marking of Section \ref{secspin}, 
$e^{2\pi i\delta''_a}=-\kappa_a\,e^{-i\pi\tau_{aa}}\,(1+O(r/d))$,
 where the left side is $\pm1$. 
 The only evidence we have is a formal one: the part of the first-order term of 
Theorem \ref{thmexp} that is odd in $\lambda_a$, namely 
$\lambda_a\bigl(\frac1{(x-w_a)(w_{-a}-y)}-\frac1{(x-w_{-a})(w_a-y)}\bigr)
=\lambda_a\frac{(x-y)(w_a-w_{-a})}{(x-w_a)(x-w_{-a})(y-w_a)(y-w_{-a})}$, 
agrees with the term of $\theta[\delta](\int_y^x\zeta)/(\theta[\delta](0)E)$ odd 
in $e^{2\pi i\delta''_a}$ to leading order, 
$e^{2\pi i\delta''_a}e^{i\pi\tau_{aa}}\,
(w_a-w_{-a})^2(x-y)/\bigl((x-w_a)(x-w_{-a})(y-w_a)(y-w_{-a})\bigr)$, 
if $\lambda_a=e^{2\pi i\delta''_a}e^{i\pi\tau_{aa}}(w_a-w_{-a})$, 
which is the claimed relation. 

\end{remark}

\section{Reduction formulas}

The theorem below describes the residue reconstruction. 
Assume the hypotheses of Theorem \ref{thmE} (in particular, (H)). 
Let $H$ be a meromorphic section of $\Lg$ on $\Sg$, with poles of arbitrary (finite) 
order at a finite set of distinct points $p_1$, $\dots$, $p_n$ of $\Sg$ 
(not of $\Omega_0$: each $p_i$ denotes a single point of the surface, 
not a choice of lift), and let $x\in\Sg\setminus\{p_1,\dots,p_n\}$.

\begin{theorem}
\label{thmred}
Then
\[
H(x)=\sum_{i=1}^n\Res_{y=p_i}\bigl(S_\eps(x,y)\,H(y)\bigr),
\]
an identity of sections of $\Lg$ in the variable $x$. For fixed $x$, 
$S_\eps(x,y)H(y)$ has ordinary scalar residues at each $p_i$.  

\end{theorem}

Note that the formula above represents a product of two sections of $\Lg$ in $y$. 
 Hence, it is a meromorphic $1$-form on $\Sg$ in $y$, by Proposition \ref{propequi}. 
  The right-hand side, as $x$ varies, assembles into a section of 
$\Lg$ because $S_\eps(x,\cdot)$ is.

\begin{proof}

By Proposition \ref{propequi} the automorphy factors $j(\tilde\mu,y)$ of $S_\eps$ and of 
$H$ cancel in the product. Thus, $\eta(y):=S_\eps(x,y)H(y)$ descends to a well-defined 
meromorphic $1$-form on $\Sg$, with poles only at $y=x$ and $y=p_i$ ($i=1$, $\dots$, $n$). 
 These are finitely many points, since $H$ has finitely many poles by hypothesis and $x$ is one 
further point. In a local coordinate near $y=x$, writing $H(y)=h(y)(dy)^{1/2}$ and 
$S_\eps(x,y)=s_\eps(x,y)(dx\,dy)^{1/2}$, 
$\eta(y)=s_\eps(x,y)h(y)\,dy\,(dx)^{1/2}=\bigl(\frac1{x-y}+O(x-y)\bigr)h(y)\,dy\,(dx)^{1/2}$ 
(Proposition \ref{propdiag} (ii)). Thus, using skew symmetry of $\frac1{x-y}$, we obtain 
$\Res_{y=x}\eta=-h(x)(dx)^{1/2}=-H(x)$. The residue theorem on the compact surface $\Sg$ 
(sum of all residues of a meromorphic $1$-form is $0$) gives $-H(x)+\sum_i\Res_{y=p_i}\eta=0$, 
i.e., the required claim.
\end{proof}

\begin{remark}
\label{remresarg}

The proof of Theorem \ref{thmred} applies the residue theorem on the compact surface $\Sg$ 
to the meromorphic $1$-form $y\mapsto S_\eps(x,y)\,\phi(y)$.  
 It is a section of $\Lg\otimes\Lg^{-1}\otimes K=K$, hence a meromorphic $1$-form for $\phi$ 
a section of $\Lg$). Its total residue is zero, and its poles are the pole $y=x$ of 
$S_\eps(x,\cdot)$ and the poles of $\phi$. There are no others since $S_\eps(x,\cdot)$ 
is holomorphic off $x$ (see Theorem \ref{thmconv}). In a chart the sum of residues is computed 
on small circles about these points. Thus, no contour at $\infty$ contributes: $\infty$ is
 an ordinary point of $\Sg$ and the homogeneous formulation (eq. \eqref{eqshat}) has no 
singularity there. Corollary \ref{corcoord} uses the same residue normalization: 
$\mathrm{Res}_{y=x_j}$ of a $1$-form is chart-independent, and the half weight factors 
transform with the same $j$-multiplier on both sides.

\end{remark}

In the coordinate form we have 

\begin{corollary}
\label{corcoord}

With $H=h(z)(dz)^{1/2}$, $c^{(i)}_m:=\Res_{y=p_i}(y-p_i)^mh(y)\,dy$ 
(thus $h=\sum_m c^{(i)}_m(y-p_i)^{-m-1}+\text{regular part}$) and $\dd_y^{(m)}=\frac1{m!}\dd_y^m$:
\[
h(x)=\sum_i\sum_{m=0}^{M_i}c^{(i)}_m\,\dd_y^{(m)}s_\eps(x,y)\big|_{y=p_i},  
\]
where $M_i$ is pole order of $H$ at $p_i$, finite by hypothesis. 
\end{corollary}

The coefficients $c_m^{(i)}$ depend on the coordinate. Theorem \ref{thmred} is the 
invariant statement. There is no holomorphic remainder because $h^0(\Lg)=0$ 
(see Theorem \ref{thmnonsing} (i)). The notation $\dd_y^{(m)}:=\frac1{m!}\dd_y^m$ 
  is for an ordinary derivative in $y$ (not an exterior/Lie derivative).

\subsection{VOSA correlation functions}

 \textit{Conventions:} $V=V_{\bar0}\oplus V_{\bar1}$ is a $\Z/2\Z$-graded VOSA, 
$|v|\in\Z/2\Z$ for $v\in V_{\overline{|v|}}$ homogeneous. A vertex operator is 
 $Y(v,x)=\sum_nv_{(n)}x^{-n-1}$. 
$\psi\in V_{1/2}\cap V_{\bar1}$ is an odd primary vector with truncation, 
$\psi_{(n)}v=0$ for $n\gg 0$  for each $v$. Locality/OPE:  
$\psi(x)Y(v,y)=\sum_{n\ge 0}(x-y)^{-n-1}Y(\psi_{(n)}v,y)+\mathrm{regular part}$,
 a finite polar sum. 
 The sum truncates by the assumption on $\psi_{(n)}v$. 
We assume also the following hypothesises. 

\smallskip\textit{Hypothesis (A).} 

For $v_1$, $\dots$, $v_n\in V$ homogeneous and  
distinct points $y_1,\dots,y_n\in\Sg$, the genus $g$ correlation function 
$H(x):=\langle\psi(x)Y(v_1,y_1)\cdots Y(v_n,y_n)\rangle_{\Sg,\eps}$ 
(i.e., the unnormalized correlation function: no division by the genus $g$ partition function 
$Z_g$ is performed, thus $Z_g\ne 0$ is not assumed;  all conditions below are linear in $H$, 
hence they hold for $H$ if and only if they hold for $H/Z_g$ whenever $Z_g\ne 0$) 
is a meromorphic section of $\Lg$ in $x$, holomorphic away from $\{y_1,\dots,y_n\}$, and 
 at each $y_k$ its principal part is the one dictated by the OPE above. In a local  
coordinate near $y_k$, we have 
\[
H(x)=\sum_{m\ge 0}(x-y_k)^{-m-1}(-1)^{\sum_{l<k}|v_l|}
\bigl\langle Y(v_1,y_1)\cdots Y(\psi_{(m)}v_k,y_k)\cdots Y(v_n,y_n)\bigr\rangle+O((x-y_k)^0). 
\]
The sign $(-1)^{\sum_{l<k}|v_l|}$ arising from anticommuting the odd field 
$\psi(x)$ past $Y(v_l,y_l)$, $l<k$ (Koszul sign rule for $\Z/2\Z$-graded vertex operators; 
no further sign occurs from the ordering of the modes $\psi_{(m)}$, 
as these act on the single vector $v_k$). 
Under condition (A), the Fermionic Zhu-type reduction is given by 

\begin{corollary}
\label{corzhu}
\[
\Bigl\langle\psi(x)\prod_kY(v_k,y_k)\Bigr\rangle
=\sum_{k=1}^n(-1)^{\sum_{l<k}|v_l|}\sum_{m\ge 0}\dd_y^{(m)}s_\eps(x,y)\big|_{y_k}
\Bigl\langle\cdots Y(\psi_{(m)}v_k,y_k)\cdots\Bigr\rangle.
\]
\end{corollary}
\begin{proof}
Theorem \ref{thmred} and Corollary \ref{corcoord}, with $c^{(k)}_m$ read from the OPE 
and the Koszul sign from moving $\psi$ past the insertions $v_l$, $l<k$.
\end{proof}

\textit{Hypothesis (A$_{\rm ff}$).} 

For the free fermion VOSA 
($\psi_{(0)}\psi=\mathbf1$, $\psi_{(m)}\psi=0$ for $m\ge 1$): for every $m\ge 0$ 
the unnormalized correlator 
$F_m(x_1,\dots,x_m)=\langle\psi(x_1)\cdots\psi(x_m)\rangle_{\Sg,\eps}$ 
is a meromorphic section of $\Lg^{\boxtimes m}$ on $\Sg^m$, holomorphic away from the 
diagonals $\{x_i=x_j\}$, with at most simple poles along the diagonals, so that in each 
variable $x_i$ it is a meromorphic section of $\Lg$ on $\Sg$, regular except at the points 
$x_j$ ($j\ne i$); and (in $x_1$) its principal part 
$\frac{(-1)^{j}}{x_1-x_j}F_{m-2}(\hat x_1,\hat x_j)$ at $x_j$; $F_0=Z_g$ is a constant, 
the genus $g$ partition function, which may be any complex number (no division by $Z_g$, 
and no assumption $Z_g\ne 0$).

\begin{proposition}
\label{proppf}

(1) Assume the hypotheses of Theorem \ref{thmE}.  
  For each constant $F_0=Z$ there is 
one family $\{F_m\}_{m\ge 0}$ satisfying (A$_{\rm ff}$), namely $F_m=0$ for $m$ odd and 
$F_{2n}=Z\cdot\mathrm{Pf}[S_\eps(x_i,x_j)]_{i<j}$.

(2) If, in addition, (A$_{\rm ff}$) holds for the correlation  
functions $F_m=\langle\psi(x_1)\cdots\psi(x_m)\rangle_{\Sg,\eps}$ of the free fermion VOSA 
 then those correlation functions equal $Z$ multiplied by 
the Pfaffians above, where $Z$ is the partition function (normalized correlators 
are obtained by dividing by $Z$ if $Z\ne 0$). 
\end{proposition}

Note that (1) is a statement purely about meromorphic sections of $\Lg$ satisfying 
the formulated analytic conditions with no reference to the genus $g$ correlation 
functions for an actual free fermion VOSA satisfying (A$_{\rm ff}$).

\begin{proof}

Suppose $\{F_m\}$ satisfies (A$_{\rm ff}$). We show by induction on $m$ that it must 
be $Z$ multiplied by the Pfaffian family,    
 and separately that this family does satisfy (A$_{\rm ff}$).  

\emph{For odd $m$:} assume $F_k=0$ for all odd $k<m$. 
 In the variable $x_1$, poles of $F_m$ can occur only at the $x_j$, 
with principal parts proportional to $F_{m-2}$, which is odd and hence $0$. Thus, $F_m$ is 
holomorphic in $x_1$, i.e., a section of $H^0(\Lg)=0$ (Theorem \ref{thmnonsing} (i)), and 
$F_m=0$.  Note that for $m=1$ there are no poles. 

\emph{For even $m$:} by Theorem \ref{thmred} in $x_1$, 
$F_m=$ $\sum_{j\ge 2}\Res_{y=x_j}S_\eps(x_1,y) F_m(y,\ldots)$ 
  $=$ $\sum_{j\ge 2}(-1)^j$ $S_\eps(x_1,x_j)$ $\cdot F_{m-2}(\hat x_1,\hat x_j)$, 
which is $Z$ multiplied by  
 the expansion of $\mathrm{Pf}[S_\eps(x_i,x_j)]$ along its first row, given 
 $F_{m-2}=Z\,\mathrm{Pf}$. For $m=2$ this gives $F_2=Z\,S_\eps(x_1,x_2)$. Conversely, 
$Z\,\mathrm{Pf}$ is skew symmetric and, 
by first-row expansion and Proposition \ref{propskew} (i), 
meromorphic with the simple poles and principal parts in every variable. 
\end{proof}

  Free fermion correlation functions are expressed 
through Szeg\H{o} kernels for specific sewing constructions of genus one and two 
in \cite{TZ1, TZ2, TZ3, TZ4}. The Szeg\H{o} kernel on a genus two surface obtained by sewing 
two tori or by self-sewing a torus. The partition function and $n$-point generating 
functions of the rank-two free fermion VOSA as determinants of such kernels, in the sewing 
parameters and with the normalization of the partition function used there.  
   Those constructions and normalizations are not the Schottky construction of this paper. 
  We do not claim here identification with those schemes.  
 Whether the Schottky-sewn correlators of the free fermion VOSA satisfy (A$_{\rm ff}$) 
is an open quesion at every genus $g\ge 2$. Proposition \ref{proppf} is therefore depending 
on condition (A$_{\rm ff}$). 

\begin{conjecture}
\label{conjzhu}

For an arbitrary VOSA $V$, hypothesis (A) holds at every genus $g\ge 2$ on (H), and 
Corollary \ref{corzhu} is the odd part of a genus $g$ Zhu recursion extending 
\cite{TW} from VOAs to VOSAs (see Remark \ref{remTW}). 
\end{conjecture}

\begin{remark}
\label{remTW}

Let us recall the Schottky recursion procedure in Tuite-Welby. 
The genus $g$ recursion in the Schottky scheme is developed by Tuite and Welby \cite{TW}, 
  for vertex operator algebras (even VOAs, with their modules). They derive a Zhu-type 
recursion for $n$-point correlation functions on a genus $g$ surface formally sewn from a 
sphere through Schottky parameters, with the Bers quasiforms of \cite{TWb} (here, the forms 
$\Psi_N$, $N\ge 2$, of integer weight) as kernels. As far as we can verify, \cite{TW} does 
not treat vertex operator superalgebras, half-integral weights, spin structures, Szeg\H{o} 
kernels or the case $N=\frac12$, and the Schottky sewing there is formal (it is not shown 
to converge in the sense of Theorem \ref{thmconv}). Thus \cite{TW} supplies the VOA model 
for Conjecture \ref{conjzhu} and for Corollary \ref{corzhu}, but the VOSA statements here 
are new and depending on conditions (A), (A$_{\rm ff}$).
\end{remark}

\section{Expansion in the Schottky multipliers}

For $p\in\mathfrak P$ set $|\rho|_{\max}:=\max_a|\rho_a|=r(p)^2$ and 
$\mathcal K_{d,p}:=\{z\in\C:|z-w_c|\ge d(p)/2\ \textit{for all}\; c\}$. 
Since $d/2\ge r_c/2$ under (H), 
$\mathcal K_{d,p}\subset\mathcal K_p$. The total space over which uniformity is claimed is 
$\mathbb G_d:=\{(p,x,y):p\in\mathfrak P,\ x,y\in\mathcal K_{d(p),p}\}$.
Expansion in $\lambda$ is given by the 

\begin{theorem}
\label{thmexp}

Assume $g\ge 2$ and (H). For every $p=(w,\lambda)\in\mathfrak P$ and all 
$x$, $y\in\mathcal K_{d,p}$, 
\[
s_\eps(x,y)=\frac1{x-y}+\sum_{a\in\Ip}\lambda_a\left[\frac{1}{(x-w_a)(w_{-a}-y)
+\rho_a}-\frac{1}{(x-w_{-a})(w_a-y)+\rho_a}\right]+R(x,y),
\]
where
\begin{equation}
\label{eqRbound}
|R(x,y)|\le\frac{64}{3}\cdot\frac{144}{143}\cdot\frac{g(2g-1)\,|\rho|_{\max}}{d^3}\ 
\le\ \frac{22\,g(2g-1)\,|\rho|_{\max}}{d^3}.
\end{equation}

\end{theorem}

The constant is absolute. The bound depends on the data only through $g$, $d$ and 
$|\rho|_{\max}$, and holds on all of $\mathbb G_d$. The bracketed terms are exact 
(given by the contribution of the words of length one). 
   $R(x,y)=O(|\lambda|_{\max}^2)$ as $|\lambda|_{\max}:=\max_a|\lambda_a|\to 0$ 
with $g$ and $d$ held fixed, so that (H) 
continues to hold along the limit. $s_\eps$ and $R$ are holomorphic in 
$(p,x,y)$ on the interior of $\mathbb G_d$ by Theorem \ref{thmconv} (i).

\begin{proof}

Length-one words are $\gt_c$, $c\in\I$. By \eqref{eqjgen} and 
$\gamma_ax-y=\frac{(w_{-a}-y)(x-w_a)+\rho_a}{x-w_a}$, the $c=a$ term is 
$\frac{\lambda_a}{(x-w_a)(w_{-a}-y)+\rho_a}$, and the $c=-a$ term 
(with $\lambda_{-a}=-\lambda_a$) is the negative of the second bracket entry. 
 Their sum is the bracket sum in the statement, and 
$R(x,y):=s_\eps(x,y)-\frac1{x-y}-[\text{bracket sum}]
=\sum_{n\ge 2}[\text{sum of the length-}n\text{ terms}]$.

\emph{Bound on $\mathcal K_{d,p}$:} let $x\in\mathcal K_{d,p}$. In the proof of Lemma 
\ref{lemword} the outermost factor satisfies 
$|j(\gt_{c_n},x)|=r_{c_n}/|x-w_{c_n}|\le 2r_{c_n}/d\le\varpi$. Now $|x-w_{c_n}|\ge d/2$
 is available, instead of the generic $\ge r_{c_n}/2$, while the remaining $n-1$ 
 factors are $\le\varpi$ as before. Hence, for every reduced word of length $n\ge 1$,
\begin{equation}\label{eqsharpj}
|j(\gt,x)|\le\varpi^n,\qquad x\in\mathcal K_{d,p}.
\end{equation}
Let now $y\in\mathcal K_{d,p}$ and $n\ge 2$. By Lemma \ref{lemword}, 
$|\gamma x-w_{-c_1}|\le r_{c_1}\varpi\le r\varpi<\frac d{8(2g-1)}\cdot\frac1{4(2g-1)}
\le\frac d{288}$ (as $g\ge 2$), while $|y-w_{-c_1}|\ge d/2$; hence 
$|\gamma x-y|\ge\frac d2-\frac d{288}=\frac{143\,d}{288}$. With \eqref{eqsharpj}, 
each term of length $n$ is bounded by $\frac{288}{143\,d}\varpi^n$, there are 
$2g(2g-1)^{n-1}$ such words, and $q:=(2g-1)\varpi\le\frac14$. Therefore, 
\[
|R(x,y)|\le\frac{288}{143\,d}\sum_{n\ge 2}2g(2g-1)^{n-1}\varpi^n
=\frac{576\,g}{143\,d\,(2g-1)}
\sum_{n\ge 2}q^n\le\frac{576\,g}{143\,d\,(2g-1)}\cdot\frac43\,q^2, 
\]
using $\sum_{n\ge 2}q^n=q^2/(1-q)\le\frac43q^2$. Since $q^2=(2g-1)^2\varpi^2
=\frac{4(2g-1)^2|\rho|_{\max}}{d^2}$,
\[
|R(x,y)|\le\frac{576\,g}{143\,d\,(2g-1)}\cdot\frac43
\cdot\frac{4(2g-1)^2|\rho|_{\max}}{d^2}=\frac{9216}{429}\cdot\frac{g(2g-1)|
\rho|_{\max}}{d^3}=\frac{64}3\cdot\frac{144}{143}\cdot\frac{g(2g-1)|\rho|_{\max}}{d^3},
\]
and $\frac{9216}{429}=21.48\ldots\le 22$. This is \eqref{eqRbound}. 
Letting $|\lambda|_{\max}\to0$ with $g$, $
d$ fixed gives $|\rho|_{\max}=|\lambda|_{\max}^2\to0$, 
hence $R=O(|\lambda|_{\max}^2)$. Holomorphicity is Theorem \ref{thmconv} (i).
\end{proof}

The bracket is skew symmetric under $x\leftrightarrow y$, as required by Proposition 
\ref{propskew} (i).
 The expansion is in the Schottky multipliers $\lambda_a$ (the multiplier of $\gamma_a$ is 
$\approx\lambda_a^2/(w_a-w_{-a})^2$). 
 It is not a sewing expansion in the vertex-algebraic 
sense. Half integer powers of a sewing parameter also occur in the sewing expansions of 
\cite{TZ1}.  $\lambda_a$ plays a role analogous to a square root of a sewing parameter.   
We do not claim identification of the two schemes. The constants of Theorem \ref{thmexp} 
are $\frac{64}3\cdot\frac{144}{143}=\frac{9216}{429}\approx21.48\le 22$.

\section{Stress tensor}
\subsection{Fermion OPE convention} 

 The free fermion VOSA is generated by $\psi\in V_{1/2}\cap V_{\bar1}$ with 
$\psi_{(0)}\psi=\mathbf1$, $\psi_{(n)}\psi=0$ ($n\ge 1$), and Virasoro vector 
$\omega_V=\half\psi_{(-2)}\psi\in V_2$. This normalization is fixed by requiring 
the standard genus zero OPE
\begin{equation}
\label{eqTpsiope}
T(z)\psi(w)\sim\frac{1/2}{(z-w)^2}\psi(w)+\frac1{z-w}\dd\psi(w),
\end{equation}
$c=\half$, which is equivalent (by the state-field correspondence and 
$\omega_V=\half\psi_{(-2)}\psi$) to 
$\psi(x)\psi(y)\sim\frac1{x-y}\mathbf1+(x-y)\bigl(\half\psi_{(-2)}\psi\bigr)(y)+O((x-y)^2)$. 
 I.e., to $:\!\dd\psi\,\psi\!:=2T$ where 
$:\!\dd\psi\,\psi\!:(y):=\lim_{x\to y}\bigl(\dd_x\psi(x)\psi(y)-\dd_x\tfrac1{x-y}\bigr)$ 
is the standard normal-ordered product. 
 We assume (A$_{\rm ff}$) at $m=2$, i.e., (by Proposition \ref{proppf}) 
$\langle\psi(x)\psi(y)\rangle_{\Sg,\eps}=Z_g\,s_\eps(x,y)(dx\,dy)^{1/2}$ for the unnormalized 
correlator. Then $\langle T\rangle=\half\langle:\!\dd\psi\,\psi\!:\rangle$, 
and when $Z_g\ne 0$ the normalized one-point function is $\langle T\rangle/Z_g$.

Finally, having (A$_{\rm ff}$) at $m=2$ and the OPE normalization 
$\psi(x)\psi(y)\sim\frac1{x-y}$, $T=\frac12{:}\partial\psi\,\psi{:}$, 
let us discuss the stress tensor. 

\begin{proposition}
\label{propT}

Under the normalization above, the unnormalized and (for $Z_g\ne 0$) normalized one-point 
functions of $T$ are
\[
\langle T(y)\rangle=Z_g\cdot\frac12\lim_{x\to y}\frac{s_\eps(x,y)-\frac1{x-y}}{x-y},
\qquad \frac{\langle T(y)\rangle}{Z_g}=\tfrac12c_1(y)=\tfrac1{24}R_\eps(y).
\]
\end{proposition}

\begin{proof}
By Proposition \ref{propdiag} (ii), $s_\eps(x,y)-\frac1{x-y}=c_1(y)(x-y)+O((x-y)^2)$ near 
$\Delta$, thus $\lim_{x\to y}\dd_x \bigl(s_\eps(x,y)-\tfrac1{x-y}\bigr)$ $=$ $c_1(y)$. I.e., 
$\langle:\!\dd\psi\,\psi\!:(y)\rangle=c_1(y)$ and, after dividing by $Z_g\ne 0$, 
$\langle T(y)\rangle/Z_g=\half c_1(y)$, which is the displayed formula. Equivalently, it is  
obtained by dividing $s_\eps-\frac1{x-y}=c_1(y)(x-y)+O((x-y)^2)$ by $(x-y)$ before taking 
$x\to y$, as in the statement. 
 Without this division by $(x-y)$ the limit would be $0$, 
thus the division is essential. By Proposition \ref{propdiag} (iii), $R_\eps=12c_1$ 
is a projective connection. Thus, $(\langle T\rangle/Z_g)^{(z)}=\half c_1^{(z)}
=w'^2\cdot\half c_1^{(w)}+\frac1{24}\{w;z\}=w'^2(\langle T\rangle/Z_g)^{(w)}
+\frac{1/2}{12}\{w;z\}$. This is the anomalous transformation law of a stress tensor 
of central charge $c=\half$ under $z\mapsto w$, consistent with \eqref{eqTpsiope}.
\end{proof}

\emph{Consistency check of sign and scalar:} with $\psi(x)\psi(y)\sim(x-y)^{-1}$ 
and $T=\frac12{:}\partial\psi\,\psi{:}$ ($c=\frac12$), the one point function is 
$\frac12\lim_{x\to y}\partial_x\bigl(\langle\psi(x)\psi(y)\rangle/Z-\frac1{x-y}\bigr)
=\frac12c_1(y)$, using $\langle\psi(x)\psi(y)\rangle=Z\,s_\eps(x,y)$ 
(hypothesis (A$_{\rm ff}$), $m=2$). Under a coordinate change the anomaly of 
$\frac12c_1=\frac1{24}R_\eps$ is $\frac1{24}\{\phi,y\}=\frac{c}{12}\{\phi,y\}$ with 
$c=\frac12$, which fixes both the sign and the scalar factor. This agrees with the 
diagonal expansion of Remark \ref{remdiagdom}.

\begin{corollary}\label{corT}
$\langle T\rangle_{\eps}/Z_g
=\tfrac1{24}\bigl(S_B+6\sum_{i,j}\dd_i\dd_j\log\theta[\delta](0)\,\zeta_i\zeta_j\bigr)$
 when $Z_g\ne 0$ (under (A$_{\rm ff}$)), by Corollary \ref{corberg}. For $g=1$ this is 
the expression of Proposition \ref{propg1} (c).
\end{corollary}

\begin{remark} 
\label{remnabla} 

Let us see that no Bers-type moduli operator exists in our formulation.
 No analogue of the Bers moduli-derivative operator is defined here: 
$S_\eps(x,w_c)$ has weight $\half$ in $w_c$ whereas $\dd_{w_c}$ has weight $-1$, 
thus a faithful analogue would require odd partners of weight $-\tfrac32$ 
(super-moduli, cf. \cite{P}). 

\end{remark}

\section{Open problems and conclusions}
 \label{secopen}

We collect here the proven and conjectural statements of this paper. 
We proved, for the $\Ng$ lifted spin structures only:  global convergence and 
meromorphicity of the half weight Poincar\'e series on (H) (estimates and gluing) and, 
qualitatively, whenever (E) holds, in any genus $g\ge1$;  periodicity and skew-symmetry; 
$H^0(\Lg)=0$; identification of (a) the series, (b) the descended section, (c) the canonical 
Szeg\H{o} kernel and the theta-function kernel; the projective connection and its relation 
to the Bergman projective connection; the genus one kernel and dictionary; residue 
reconstruction for arbitrary pole orders; the expansion in the Schottky multipliers with 
remainder; the stress-tensor formula.

   Beyond (E): for $\delta(\Gamma)<\half$ convergence is proven 
(Theorem \ref{thmE}). It is  an open quesion whether (E) fails when $\delta(\Gamma)\ge\half$, 
and whether a Hecke-type regularization $\sum_{\gt}|j(\gt,x)|^{2s}(\cdots)$, continued 
analytically to $s=0$, extends $S_\eps$ across the locus where (E) fails, for instance,    
 across the theta-null divisors $\theta[\delta(\eps)](0)=0$.
 The dictionary for $g\ge 2$: determine $\eps\mapsto\delta''(\eps)$. Remark  
\ref{remdict} gives the expected form. 

 Grandint Hypothesis (A$_{\rm ff}$) for the Schottky scheme:  whether the Schottky-sewn  
correlation functions of the free fermion VOSA satisfy (A$_{\rm ff}$) is (to our knowledge)
 an open question at every genus $g\ge 2$. The genus one and genus two sewing schemes of 
\cite{TZ1, TZ2, TZ3, TZ4} are different constructions, and no identification of the schemes is 
claimed here.  Corollary \ref{corzhu} dependes on (A) and 
Proposition \ref{proppf} on (A$_{\rm ff}$), neither of which is established for 
Schottky-sewn correlators. 

 The other $2^{2g}-\Ng$ spin structures: 
(both even structures outside the $\Ng$-family and all odd ones), zero modes, super-moduli.
  We do not treated it in this paper. 
 It presumably requires either a different choice of generators/markings 
for each spin structure, or the super-Schottky formalism of \cite{P}. 
Super-moduli analogue of the Bers operator: see Remark \ref{remnabla}. 

It remains conjectural: the dictionary for $g\ge 2$, continuation 
beyond (E), and the general-VOSA recursion (Conjecture \ref{conjzhu}).  
 
The material of this paper is also useful in other areas of condensed matter theory 
\cite{chernodub2017scale, zhang2020influence, zubkov2012momentum, zubkov2017topology, kmmzz}, 
and foliation theory \cite{BG, BGG, BEG, Z1, Z2, Z3}. 

\section*{Acknowledgments}

The author is supported by the Institute of Mathematics, Academy of Sciences of the Czech
Republic (RVO 67985840). 

	\medskip
	\noindent\textbf{Data Availability.}
	Data sharing is not applicable to this article as no datasets were generated
	or analysed during the current study.

	\medskip
	\noindent\textbf{Declarations}

This paper arose in close collaboration with ChatGPT, an AI assistant developed by 
OpenAI. The collaboration involved exploring proof strategies, refining arguments, locating
references, and drafting. Revision and writing of final text have been done by the author.
The author takes full responsibility for the mathematical content and the final text.

\medskip
	\noindent\textbf{Conflict of interest.}
	The authors have no conflicts of interest to declare that are relevant to the
	content of this article.


\end{document}